\documentclass[a4paper]{article}

\usepackage{arxiv}

\usepackage{bookmark}
\usepackage{multirow,multicol,diagbox,makecell}
\usepackage{booktabs}
\usepackage{graphicx}
\usepackage{subfigure}
\usepackage{amsmath,amsthm}
\usepackage{mathtools,esint,physics,amssymb}
\usepackage{amsfonts}
\usepackage{enumerate}
\usepackage{url}
\usepackage{xcolor}
\usepackage{tikz}
\usepackage{siunitx}
\usepackage[capitalise,nameinlink,noabbrev]{cleveref}
\usepackage{hyperref}
\usetikzlibrary{graphs, quotes, positioning, shapes.geometric, arrows.meta, fit, positioning, shapes.geometric, shapes.arrows}

\numberwithin{equation}{section}

\newcommand{\R}{\mathbb{R}}
\newcommand{\uns}{\underline{n}^2} % underline n square
\newcommand{\uno}{\underline{n}_0} % underline n sub 0
\newcommand{\bp}{\mathbb{P}}
\newcommand{\be}{\mathbb{E}}
\DeclareMathOperator{\dist}{dist}

\newtheorem{theorem}{Theorem}
\newtheorem{remark}{Remark}
\newtheorem{lemma}{Lemma}
\newtheorem{definition}{Definition}
\newtheorem{corollary}{Corollary}

\title{Sediment Concentration Estimation via Multiscale Inverse Problem and Stochastic Homogenization}

\renewcommand{\shorttitle}{Sediment Measurement via Multiscale Inverse Problem}

\author{Jiwei Li$^{1}$, Lingyun Qiu$^{2,3}$, Zhongjing Wang$^{4,5}$, Hui Yu$^{6}$, Siqin Zheng$^{7}$\\
$^1$ Department of Mathematical Sciences, Tsinghua University, Beijing 100084, China\\
(\href{mailto:li-jw20@mails.tsinghua.edu.cn}{li-jw20@mails.tsinghua.edu.cn})\\
$^2$ Yau Mathematical Sciences Center, Tsinghua University, Beijing 100084, China\\
$^3$ Yanqi Lake Beijing Institute of Mathematical Sciences and Applications, Beijing 101408, China\\
(\href{mailto:lyqiu@tsinghua.edu.cn}{lyqiu@tsinghua.edu.cn})\\
$^4$ Department of Hydraulic Engineering, Tsinghua University, Beijing 100084, China\\
$^5$ Breeding Base for State Key Laboratory of Land Degradation and Ecological Restoration in Northwest China, \\
Ningxia University, Yinchuan 750021, Ningxia, China \\
(\href{mailto:zj.wang@tsinghua.edu.cn}{zj.wang@tsinghua.edu.cn})\\
$^6$ School of Mathematics and Computational Science, Xiangtan University, Xiangtan 411105, Hunan, China\\
(\href{mailto:huiyu@xtu.edu.cn}{huiyu@xtu.edu.cn})\\
$^7$ Department of Mathematics, City University of Hong Kong, Kowloon, Hong Kong SAR, China\\
(\href{mailto:siqzheng@cityu.edu.hk}{siqzheng@cityu.edu.hk})
}

\begin{document}
\maketitle

\begin{abstract}
	We develop a multiscale framework for estimating sediment concentration in water flow from acoustic wave measurements. At the microscopic scale, the sediment distribution is modeled by a spatially inhomogeneous Poisson cloud, while the quantity of interest is its macroscopic concentration. For the associated random wave model, we derive an effective medium whose coefficient is explicitly related to the local probability of sediment occurrence. This effective description avoids resolving individual sediment particles and provides a computationally tractable forward model for inversion. We then formulate the recovery of the effective medium, and hence the sediment concentration, as an inverse medium problem from partial boundary measurements, and investigate numerical strategies including model mollification and shot averaging. Numerical experiments demonstrate that the effective model captures the macroscopic wave behavior and can be used to obtain accurate estimates of sediment concentration.
\end{abstract}

\keywords{
    sediment concentration measurement, multiscale modeling, inverse problem, homogenization theory
}

\noindent{\bfseries \emph{MSC Classification}}\enspace{35B27, 35R30, 76M50}

\section{Introduction}\label{sec_intro}
Sediment measurement stands as a cornerstone in environmental studies, hydraulic engineering, and geoscience research, as it provides crucial insights into the dynamics and characteristics of aquatic systems. Accurate quantification of sediment concentration in water bodies is essential for understanding sediment transport processes, assessing water quality, and designing effective erosion control measures. However, the inherent complexity of sediment dynamics coupled with limitations in measurement techniques poses significant challenges to achieving precise and reliable sediment concentration estimates.

Historically, sediment concentration measurement techniques have predominantly relied on direct sampling methods \cite{grayComparabilityAccuracyFluvialsediment2002,thorneReviewAcousticMeasurement2002,wrenFieldTechniquesSuspendedSediment2000} such as grab sampling, core sampling, and filtration. While these methods offer tangible measurements, they are often labor-intensive, time-consuming, and susceptible to spatial and temporal variability, limiting their applicability in dynamic environments. Moreover, these techniques may disrupt the natural sediment distribution and fail to capture fine-scale variations, leading to potential inaccuracies in concentration estimates.

In response to the limitations of direct sampling methods, a variety of indirect measurement approaches have been developed, including acoustic Doppler techniques \cite{meralStudyEstimatingSediment2015,crawfordDeterminingSuspendedSand1993,thorneAcousticMeasurementsSuspended1997,thostesonSimplifiedMethodDetermining1998a,pedocchiAcousticMeasurementSuspended2012,sEstimationSuspendedSediment2009,toppingLongtermContinuousAcoustical2016a},
optical methods \cite{schoellhamerContinuousMonitoringSuspended2002,greenMeasurementConstituentConcentrations1993,blackSuspendedSandMeasurements1994,wrenFieldTechniquesSuspendedSediment2000},
and remote sensing technologies \cite{wrenFieldTechniquesSuspendedSediment2000,novoEffectSedimentType1989,novoEffectViewingGeometry1989,choubeyEffectPropertiesSediment1994,gaoRoleSpatialResolution1997}.
These techniques infer sediment concentration through its influence on acoustic or optical signals and therefore allow non-invasive measurements over considerably larger spatial and temporal scales. Nevertheless, relating the measured signals quantitatively to sediment concentration remains challenging, particularly when the sediment distribution exhibits strong microscopic heterogeneity.

In this work, we consider the sediment-laden flow from a multiscale acoustic perspective. Individual sediment particles generate rapid spatial variations of the acoustic properties at a small characteristic scale, whereas the quantity relevant for sediment monitoring is the macroscopic concentration. Resolving every particle in either forward simulation or inversion is computationally prohibitive and, more importantly, unnecessary if only the macroscopic sediment distribution is sought. This motivates replacing the microscopic random medium by an effective medium that preserves its large-scale influence on wave propagation.

Homogenization provides a systematic framework for such an effective description and has been extensively studied for wave propagation in heterogeneous media. For example, homogenized models for seismic and elastic wave equations in heterogeneous media have been investigated \cite{capdeville1DNonperiodicHomogenization2010,capdeville2DNonperiodicHomogenization2010,guillot2DNonperiodicHomogenization2010a},
together with higher-order extensions in related settings \cite{capdevilleSecondOrderHomogenization2007,chenDispersiveModelWave2000,fishHigherOrderHomogenizationInitial2001,fishNonlocalDispersiveModel2002a}.
Random Helmholtz equations and the fluctuations around their homogenized limits have also been studied using corrector theory; see, for instance, \cite{Jing2019}. These works establish the broader homogenization framework on which the present study is built. Our purpose here is not to introduce a new general homogenization mechanism, but rather to develop and quantitatively justify an effective description tailored to a stochastic model of sediment-laden flow.

Effective models have likewise played an important role in multiscale inverse problems. General strategies for parameter estimation when the observations contain fine-scale features but the unknown of interest is macroscopic have been considered in \cite{eHeterognousMultiscaleMethods2003,frederickNumericalMethodsMultiscale2017,nolenMultiscaleModellingInverse2012}. The influence of unresolved random fluctuations on parameter estimation was analyzed in \cite{nolenFineScaleUncertainty2009}, while a Bayesian numerical homogenization approach for the recovery of macroscopic parameters from multiscale observations was developed in \cite{abdulleBayesianNumericalHomogenization2020}. These studies show that replacing a fine-scale forward model by an appropriate effective model can substantially reduce the computational complexity of inversion, while the discrepancy between the heterogeneous data and the effective forward map must be taken into account.

Against this background, the present work focuses on a specific modeling question arising in sediment concentration measurement: how can such a macroscopic quantity be encoded by a microscopic random particle distribution, and how can it be recovered from acoustic measurements without resolving the individual particles? To address these questions, we model the centers of the sediment particles by a spatially inhomogeneous Poisson point process with intensity function $\rho(x)$. Particles are represented by balls of characteristic radius $\epsilon$. This construction yields an explicit local probability $p(x) = 1-e^{-\rho(x)}$ for the presence of sediment in the small-scale limit. The resulting acoustic coefficient is therefore random and rapidly varying at the microscopic level, while $p(x)$ provides a deterministic macroscopic description directly related to sediment concentration.

The contributions of this work are threefold. First, for this spatially inhomogeneous Poisson-cloud model, we derive the corresponding effective acoustic coefficient and establish a quantitative estimate for the approximation of the heterogeneous wave field by the effective field. In the frequency domain, the randomness appears in the zeroth-order refractive-index term, so the effective coefficient can be identified explicitly through its local statistical average; the role of the homogenization analysis is therefore to justify this effective replacement quantitatively for the nonstationary particle model considered here. Second, the effective slowness square admits an explicit relation with $p(x)$, allowing the macroscopic sediment distribution, and consequently the sediment concentration, to be obtained once the effective acoustic medium has been reconstructed. Third, we formulate this reconstruction as an inverse problem using partial boundary wave measurements and investigate numerical strategies adapted to the multiscale data, including model mollification and averaging over different realizations of the microscopic sediment configuration. Numerical experiments are used to assess both the effective approximation and the resulting sediment concentration estimates.

The remainder of the paper is structured as follows. \Cref{sec_review} reviews several conventional techniques for sediment concentration measurement. In \cref{sec_model}, we introduce the spatially inhomogeneous Poisson-cloud model, derive and quantitatively justify its effective acoustic medium, and formulate the associated multiscale inverse problem. \Cref{sec_numerical} outlines the numerical methods proposed to solve the multiscale inverse problem. Numerical experiments validating the effective approximation and the sediment concentration reconstruction are presented in \cref{sec_results}. \Cref{sec_conclusions} concludes the paper and discusses possible extensions.

\section{Review of sediment concentration measurement methods}\label{sec_review}
In this section, we provide a brief overview of several conventional instruments and associated technologies used to measure the sediment concentration. These methods include direct sampling \cite{grayComparabilityAccuracyFluvialsediment2002,wrenFieldTechniquesSuspendedSediment2000,thorneReviewAcousticMeasurement2002},
acoustic Doppler techniques \cite{meralStudyEstimatingSediment2015,crawfordDeterminingSuspendedSand1993,thorneAcousticMeasurementsSuspended1997,thostesonSimplifiedMethodDetermining1998a,pedocchiAcousticMeasurementSuspended2012,sEstimationSuspendedSediment2009,toppingLongtermContinuousAcoustical2016a},
optical methods \cite{schoellhamerContinuousMonitoringSuspended2002,greenMeasurementConstituentConcentrations1993,blackSuspendedSandMeasurements1994,wrenFieldTechniquesSuspendedSediment2000},
and remote sensing technologies \cite{wrenFieldTechniquesSuspendedSediment2000,novoEffectSedimentType1989,novoEffectViewingGeometry1989,choubeyEffectPropertiesSediment1994,gaoRoleSpatialResolution1997}.
In the following, we will discuss the principles, advantages, and limitations of each sort of methods, highlighting the challenges and opportunities in sediment concentration measurement.

Direct sampling methods are among the oldest and most commonly used techniques for measuring the sediment concentration in water. These methods involve collecting sediment samples using tools such as sediment samplers, drying and weighing the sediment to calculate the volume concentration as $m/(V\rho_s)$, where $m$ is the mass of sediment collected, $V$ is the volume of sampled water and $\rho_s$ is the density of sediment. While providing direct measurements, these methods are labor-intensive, time-consuming, and limited in spatial and temporal coverage. Moreover, their sensitivity is influenced by sampling and laboratory procedures \cite{wrenFieldTechniquesSuspendedSediment2000,schoellhamerContinuousMonitoringSuspended2002}.

Acoustic Doppler meters, originally designed for flow-velocity measurements, can also estimate suspended sediment concentration by analyzing the intensity of backscattered sound waves. The relationship between backscatter intensity $I$ and sediment concentration $C$ is expressed as $I=I_0e^{-\alpha C}$, where $I_0$ is the intensity of incident sound, and $\alpha$ is the attenuation coefficient determined through calibration (see details in \cite{hanesAcousticMeasurementsSuspended1988,pomaziAcousticBasedAssessment2022}).
Researchers have also explored using signal-to-noise ratios from acoustic Doppler velocimeters to estimate sediment concentration \cite{sEstimationSuspendedSediment2009,hosseiniSynchronousMeasurementsVelocity2006a,salehiUsingVelocimeterSignal2011,haUsingADVBackscatter2009}.
Also, the framework used to interpret the backscattered signals from suspensions of non-cohesive sediments can be found in \cite{thorneReviewAcousticMeasurement2002}.
While offering non-invasive, continuous measurements over large areas and depths, these techniques face challenges in translating backscatter data into sediment concentrations and may struggle to capture fine-scale variations.

Optical methods measure sediment concentration based on light attenuation, using sensors to detect light absorbed or scattered by suspended particles. Two common techniques, optical backscatter (OBS) and optical transmission, rely on infrared or visible light. OBS measures the backscattered light detected by photodiodes positioned around the emitter, with backscatter strength used to estimate the sediment concentration. Studies show a nearly linear response of OBS to homogeneous sediment concentrations over a wide range \cite{greenMeasurementConstituentConcentrations1993,blackSuspendedSandMeasurements1994}.
Optical transmission measures light attenuation as it passes through the sample, offering higher sensitivity to low particle concentrations. However, both methods are affected by the particle size \cite{blackSuspendedSandMeasurements1994,xuConvertingNearbottomOBS1997,ludwigLaboratoryEvaluationOptical1990,cliffordLaboratoryFieldAssessment1995}.

Remote sensing technologies, including satellite and aerial imagery, estimate sediment concentration in large water bodies based on spectral reflectance, where reflected light is related to water properties like suspended particle concentration. Studies have explored the correlation between sediment concentration and reflected radiation \cite{novoEffectViewingGeometry1989,novoEffectSedimentType1989,bhargavaLightPenetrationDepth1991,choubeyEffectPropertiesSediment1994}.
While these methods provide spatially continuous measurements over large areas, they may struggle to capture small variations in sediment concentration. Additionally, large errors can occur when the sediment type is unknown due to the strong correlation between spectral readings and sediment mineral composition \cite{novoEffectSedimentType1989,choubeyEffectPropertiesSediment1994}.

\section{Mathematical model}\label{sec_model}
In this section, we present the mathematical model for the sediment concentration measurement in the water flow. We first introduce a stochastic microscopic model for the sediment distribution and the resulting acoustic wave equation. We then derive and quantitatively justify an effective medium model for this heterogeneous system. Finally, we formulate the recovery of the macroscopic sediment distribution from boundary wave measurements as a multiscale inverse problem.

As in acoustic Doppler techniques, transmitters and receivers are placed on the boundary of the measurement domain, and the sediment concentration is inferred through its influence on acoustic wave propagation. The distinction here is that, rather than relying on a direct empirical calibration between signal attenuation and sediment concentration, we model the heterogeneous wave propagation explicitly and recover the macroscopic sediment distribution through an inverse medium problem.

\subsection{Acoustic structure of sediment-laden flow}\label{sec_physical_model_sediment}
Firstly, we introduce a physical model to describe the sediment concentration in water flow using the p-wave velocity. Drawing inspiration from acoustic Doppler-based instruments and the distinct acoustic properties of sediment and water, we model the sediment distribution as a highly heterogeneous acoustic medium, denoted by $c_\epsilon(x,t)$. This function represents the p-wave velocity, which varies depending on the medium: at a given location $x$ and time $t$, $c_\epsilon(x,t)$ corresponds to the p-wave velocity in sediment if sediment is present; otherwise, it reflects the p-wave velocity in water.

Next, consider the sediment distribution in a river as an example. At the macroscopic scale, sediment particles are primarily concentrated near the riverbed, with the concentration gradually decreasing as the depth decreases. At the microscopic scale, however, individual sediment particles are distributed randomly and exhibit rapid variations. This variability highlights the nature of the p-wave velocity $c_\epsilon(x,t)$, which reflects the heterogeneous distribution of sediment as a multiscale field.

To describe this multiscale behavior, we introduce a smooth function $p(x,t)$ to represent the probability of sediment present at location $x$ and time $t$. The instantaneous sediment concentration can then be expressed as the integral of $p(x,t)$ with respect to $x$. Notably, the propagation speed of the p-wave is much faster than the sediment transport speed, allowing us to approximate the sediment distribution as a static field during wave propagation. Thus, we have $p(x,t)\approx p(x)$ and $c_\epsilon(x,t)\approx c_\epsilon(x)$.

More formally, the heterogeneous p-wave velocity is modeled as the Poisson cloud model. It is based on a heterogeneous Poisson point process, a widely used random process that describes the random distribution of particles in the space \cite{Kingman1993}.
\begin{definition}\label{def_cloud}
    A Poisson point process on $\R^d$ with a heterogeneity scale $\epsilon>0$ and a nonnegative intensity function $\rho$ is a random countable subset $\Pi_\epsilon$ of $\R^d$, satisfying that
    \begin{enumerate}
        \item for any Borel subset $A$ of $\R^d$, the number of points in $\Pi_\epsilon\cap A$, denoted by $N_\epsilon(A)$, is a random variable following the Poisson law with mean
        \begin{equation*}
            \lambda_\epsilon(A) = \frac{1}{V_d \epsilon^d}\int_A \rho(x)\,dx,
        \end{equation*}
        where $V_d$ is the volume of the unit ball in $\R^d$, and 
        \item for any disjoint Borel subsets $A_1,\cdots,A_k$ of $\R^d$, the random variables $N_\epsilon(A_1),\cdots$, $N_\epsilon(A_k)$ are independent.
    \end{enumerate}
\end{definition}

The parameter $\epsilon$ is the characteristic length scale, which represents the size of the small-scale sediments. Intuitively, the smaller the sediment size $\epsilon$, the more sediments are in the same set $A$ on average, so it is reasonable in the definition to divide the intensity function $\rho$ by the volume of a ball with radius $\epsilon$.

Viewing the sediments as spherical particles of radius $\epsilon$, denoted by $B(x,\epsilon)$ with $x$ representing their locations, we define the random field of the heterogeneous p-wave velocity over the space as
\begin{equation}\label{eq_c_eps}
    c_\epsilon(x) = c_\epsilon(x,\Pi_\epsilon) = \begin{cases}
        c_1, & \operatorname{dist}(x,\Pi_\epsilon) < \epsilon, \\
        c_0, & \operatorname{dist}(x,\Pi_\epsilon)\ge\epsilon,
    \end{cases}
\end{equation}
where $c_0$ and $c_1$ denote the p-wave velocity of water and sediments respectively, which are assumed to be constant for simplicity. Therefore, the probability of sediments arising at location $x\in\R^d$ is
\begin{equation}
    p_\epsilon(x) = 1 - \bp[N_\epsilon(B(x,\epsilon)) = 0] = 1 - e^{-\fint_{B(x,\epsilon)} \rho(y)\,dy}.
\end{equation}
Since $\fint_{B(x,\epsilon)} \rho(y)\,dy \to \rho(x)$ as $\epsilon\to 0$ for continuous $\rho$, it is natural to define
\begin{equation}
    p(x) = \lim_{\epsilon\to 0} p_\epsilon(x) = 1 - e^{-\rho(x)}
\end{equation}
as the probability of sediments arising at location $x\in\R^d$ at the macroscopic scale.

\subsection{Wave propagation in heterogeneous medium}
Consider the highly heterogeneous velocity profile $c_\epsilon(x)$ representing the p-wave velocity at $x$ in the flow. For simplicity, we assume that the velocity profile $c_\epsilon$ varies rapidly in a bounded domain $D\subset\mathbb{R}^d(d=2,3)$ and is constant outside $D$. It means that the intensity function $\rho(x)$ defining the Poisson cloud model vanishes outside $D$, i.e.\ $\operatorname{supp}\rho\subset\bar{D}$. Given a known acoustic source $f(x,t)$, the equation governing the wave propagation is given by the following wave equation:
\begin{subequations}\label{eq_heterogeneous_wave}
    \begin{align}
        \left(\frac{1}{c_\epsilon^2(x)}\partial_t^2-\Delta\right)u_\epsilon=f, &\quad x\in \mathbb{R}^d,\ t\in(0,T),\\
        u_\epsilon(x,0)=0, &\quad x\in \mathbb{R}^d,\\
        \partial_t u_\epsilon(x,0)=0, &\quad x\in \mathbb{R}^d.
    \end{align}
\end{subequations} 
Let the slowness square $m_\epsilon(x)$ be $1/c_\epsilon^2(x)$. The solution $u_\epsilon(x,t)$ of \eqref{eq_heterogeneous_wave} represents the wave field propagating in the heterogeneous medium. It is measured by the receivers located at the partial boundary $\Gamma\subset\partial D$, and the measurements $u_\epsilon|_{\Gamma\times(0,T)}$ are used to infer the sediment distribution $p(x)$. For convenience, we rewrite the measurement of the wave field as an operator 
\begin{equation}
    {\cal F}_{\rm het}^\epsilon\colon m_\epsilon(x)\mapsto u_\epsilon|_{\Gamma\times(0,T)}.
\end{equation}

However, the highly heterogeneous nature of the p-wave velocity $c_\epsilon(x)$ makes the direct solution of \eqref{eq_heterogeneous_wave} computationally expensive, which means that a quite fine mesh is required to capture the small-scale heterogeneity. Moreover, we are interested in the macroscopic sediment distribution $p(x)$ instead of any concrete microscopic structure of $c_\epsilon(x)$. To avoid recovering specific realizations of $c_\epsilon(x)$, we consider the stochastic homogenization theory to derive the effective medium model, where the effective coefficient directly matches the sediment distribution $p(x)$. In the following, we will present the stochastic homogenization theory based on the wave equation \eqref{eq_heterogeneous_wave} in the frequency domain, i.e.\ the Helmholtz equation.

\subsection{Stochastic homogenization of wave equations}\label{sec_homogenization}
Suppose that the acoustic source $f$ is in the form of
\begin{equation*}
    f(x,t) = \int_I f_\omega(x) e^{i\omega t}\,d\omega
\end{equation*}
with a bounded interval $I\subset\subset(0,\infty)$. This assumption is quite practical, since the acoustic source we use in practice always has a finite frequency band. By linearity, the solution to the wave equation \eqref{eq_heterogeneous_wave} can be decomposed as
\begin{equation*}
    u_\epsilon(x,t) = \int_I u_{\epsilon,\omega}(x) e^{i\omega t}\,d\omega,
\end{equation*}
with the components $u_{\epsilon,\omega}$ solving the Helmholtz equation
\begin{equation}\label{eq_hete_Helmholtz}
    \left(-\omega^2 m_\epsilon(x)-\Delta\right)u_{\epsilon,\omega} = f_\omega, \quad x\in\R^d.
\end{equation}
For simplicity in the following analysis, we generalize the coefficients of \eqref{eq_hete_Helmholtz} to the complex domain, omit the subscript $\omega$ and rewrite it as
\begin{equation}\label{eq:hete}
    -\left(\Delta + k^2\uns_\epsilon(x)\right)u_\epsilon = f, \quad x\in\R^d,
\end{equation}
where $k=\omega/c_0$ is the wave number with frequency $\omega$ in water, and $\uns_\epsilon(x)$ is the complex refractive index, which equals $\uns_0 = 1 + i\kappa_0$ in water and $\uns_1 = n_1^2 + i\kappa_1$ in sand. Here $n_1 = c_0/c_1$ denotes the real refractive index between water and sand, and we add the absorption $\kappa_0,\kappa_1>0$ in water and sand, to bypass the analysis involving the radiation condition at infinity. As a result, any source $f\in L^2(\R^d)$ immediately gives rise to a unique solution to \eqref{eq:hete} in $L^2(\R^d)$. For the sake of simplicity, we have omitted the subscript $\omega$ and follow the same notation $u_\epsilon$ and $f$ in \eqref{eq_heterogeneous_wave}, which would not lead to any confusion since $u_\epsilon$ and $f$ only depend on $x$ in the analysis of \cref{sec_homogenization} and \cref{sec_convergence}.

We next derive an effective description of the random Helmholtz equation associated with the Poisson-cloud model. We emphasize that the random refractive index in \eqref{eq:hete} appears as a zeroth-order coefficient, while the principal part remains the deterministic Laplacian. Consequently, unlike homogenization problems with rapidly oscillating principal coefficients, no cell problem is needed to identify the leading-order effective coefficient: it is determined explicitly by the local statistics of the random medium. The purpose of the analysis below is instead to quantitatively justify this effective replacement for the spatially inhomogeneous Poisson cloud introduced in \cref{def_cloud}. Our argument is related to the corrector approach for random Helmholtz equations in \cite{Jing2019}, but exploits the finite-range dependence inherited from the Poisson construction to obtain an explicit convergence estimate.
More precisely, consider the random coefficient field
\begin{equation*}
    \uns_\epsilon(x) = \uns_\epsilon(x,\Pi_\epsilon) = \begin{cases}
        \uns_1, & \operatorname{dist}(x,\Pi_\epsilon) < \epsilon, \\
        \uns_0, & \operatorname{dist}(x,\Pi_\epsilon)\ge\epsilon.
    \end{cases}
\end{equation*}
The quantitative effective-medium result is stated as follows.
\begin{theorem}\label{thm_homogenization}
    Suppose that the nonnegative intensity function $\rho\in C^\alpha(\R^d)$ for some $\alpha\in(0,1]$ with the $\alpha$-H\"older seminorm
    \begin{equation*}
        [\rho]_\alpha = \sup_{x\ne y} \frac{|\rho(x) - \rho(y)|}{|x-y|^\alpha}.
    \end{equation*}
    Define
    \begin{align}
        \mu_\epsilon(x) &\coloneqq \be[\uns_\epsilon] = p_\epsilon(x)\uns_1 + (1-p_\epsilon(x))\uns_0, \label{eq:mu_eps} \\
        \mu(x) &\coloneqq \lim_{\epsilon\to 0} \mu_\epsilon(x) = p(x)\uns_1 + (1-p(x))\uns_0. \label{eq:mu}
    \end{align}
    For $f\in L^2(\R^d)$, let $\tilde{u}_\epsilon$ and $u$ solve
    \begin{align}
        -\left(\Delta + k^2 \mu_\epsilon\right)\tilde{u}_\epsilon &= f \ \text{in}\ \R^d, \label{eq:inter} \\
        -\left(\Delta + k^2 \mu\right)u &= f \ \text{in}\ \R^d. \label{eq:homo}
    \end{align}
    Then for any $\epsilon\in(0,1)$,
    \begin{align}
        \be\|u_\epsilon - \tilde{u}_\epsilon\|_{L^2(\R^d)}^2 &\le C(k,\uns_0,\uns_1,d) \epsilon^d \|f\|_{L^2(\R^d)}^2, \label{eq:est1} \\
        \|\tilde{u}_\epsilon - u\|_{L^2(\R^d)} &\le C(k,\uns_0,\uns_1,d,\alpha,[\rho]_\alpha) \epsilon^\alpha \|f\|_{L^2(\R^d)}, \label{eq:est2}
    \end{align}
    and consequently,
    \begin{equation}\label{eq:homo_est}
        \be\|u_\epsilon - u\|_{L^2(\R^d)}^2 \le C(k,\uns_0,\uns_1,d,\alpha,[\rho]_\alpha) \epsilon^{2\alpha}\|f\|_{L^2(\R^d)}^2.
    \end{equation}
\end{theorem}

The estimates above distinguish two different sources of approximation error. The difference $u_\epsilon-\tilde{u}_\epsilon$ describes the realization-dependent fluctuation of the heterogeneous wave field, whereas $\tilde{u}_\epsilon-u$ is the deterministic error caused by replacing the microscopic probability $p_\epsilon$ by its macroscopic limit $p$. Hence the intermediate field $\tilde{u}_\epsilon$ separates random fluctuations from the deterministic scale bias.

\cref{thm_homogenization} provides a quantitative justification for replacing the heterogeneous frequency-domain model \eqref{eq_hete_Helmholtz} by the effective equation
\begin{equation}\label{eq_eff_Helmholtz}
    \left(-\omega^2 m(x)-\Delta\right)u_\omega = f_\omega, \quad x\in\R^d,
\end{equation}
where
\begin{equation}\label{eq_m_eff}
    m(x) = \lim_{\epsilon\to 0} \mathbb{E}[m_\epsilon(x)] = \frac{p(x)}{c_1^2} + \frac{1-p(x)}{c_0^2}
\end{equation}
is the effective slowness square, which depends explicitly on the macroscopic probability $p(x)$ of sediment occurrence. Since the source $f$ is restricted to a bounded frequency band $I\subset\subset(0,\infty)$, superposition of the effective frequency components
\begin{equation*}
    u(x,t)\coloneqq \int_I u_\omega(x) e^{i\omega t}\,d\omega
\end{equation*}
leads naturally to the time-domain effective model
\begin{subequations}\label{eq_effective_wave}
    \begin{align}
        \left(m(x)\partial_t^2-\Delta\right)u=f, &\quad x\in \mathbb{R}^d,\ t\in(0,T),\\
        u(x,0)=0, &\quad x\in \mathbb{R}^d,\\
        \partial_t u(x,0)=0, &\quad x\in \mathbb{R}^d.
    \end{align}
\end{subequations}
This model contains no microscopic particle-scale oscillations and is therefore substantially less expensive to solve numerically than the heterogeneous model \eqref{eq_heterogeneous_wave}. More importantly for the inverse problem, its coefficient $m(x)$ is directly related to the macroscopic quantity $p(x)$ that we aim to reconstruct.

For later use in the inverse problem, we denote the forward operator associated with the effective model by
\begin{equation}
    {\cal F}_{\rm eff}\colon m(x)\mapsto u|_{\Gamma\times(0,T)}.
\end{equation}
Since the inverse problem is formulated in terms of measurements on $\Gamma$, while the convergence in \cref{thm_homogenization} is stated in the whole-space $L^2$ norm, we next present the corresponding convergence result for the boundary data.
\begin{corollary}\label{cor:homo}
    Let $D\subset\mathbb R^d$ be a bounded Lipschitz domain and $\Gamma\subset\partial D$ be the measurement boundary. Then under the assumptions of \cref{thm_homogenization}, for any $\epsilon\in(0,1)$,
    \begin{align}
        \be\|u_\epsilon - \tilde{u}_\epsilon\|_{H^1(\R^d)}^2 &\le C(k,\uns_0,\uns_1,d) \epsilon^{d/2} \|f\|_{L^2(\R^d)}^2, \label{eq:est1_H1} \\
        \|\tilde{u}_\epsilon - u\|_{H^1(\R^d)} &\le C(k,\uns_0,\uns_1,d,\alpha,[\rho]_\alpha) \epsilon^\alpha \|f\|_{L^2(\R^d)}, \label{eq:est2_H1}
    \end{align}
    and consequently,
    \begin{equation}\label{eq:homo_est_bdry}
        \be\|u_\epsilon - u\|_{L^2(\Gamma)}^2 \le C(k,\uns_0,\uns_1,d,\alpha,[\rho]_\alpha,D) \epsilon^{\min\{d/2, 2\alpha\}}\|f\|_{L^2(\R^d)}^2.
    \end{equation}
\end{corollary}

This result transfers the error decomposition in \cref{thm_homogenization} to the boundary measurements. It is important for the averaging strategy introduced in \cref{sec_solveIP}: averaging independent realizations can reduce the random fluctuation error, although it cannot decrease the deterministic scale error.

\cref{cor:homo} provides a quantitative justification for approximating the heterogeneous boundary measurements by those generated by the effective model. To numerically assess this approximation, in \cref{sec_results} we generate the heterogeneous and effective media on different spatial scales and compare the corresponding wave fields. The numerical results show that the effective model captures the dominant macroscopic features of wave propagation in the heterogeneous medium, thereby complementing the theoretical boundary-data estimate above.

\subsection{Proofs of the quantitative homogenization results}\label{sec_convergence}
For ease of notations, we shall abbreviate the function spaces $L^p(\R^d)$ and $H^s(\R^d)$ with $L^p$ and $H^s$ respectively, and use $C(\cdot)$ to denote the constant with dependence in parentheses. The complex conjugate of a complex number $z$ is denoted by $z^*$. Before proving the homogenization results, we prepare with two lemmas.
\begin{lemma}
    Let $k>0$ and $q\in L^\infty$ with
    \begin{equation*}
        \inf_{x\in\R^d} \Im q(x) = c > 0.
    \end{equation*}
    For $g\in L^2$, let $v\in L^2$ be a distributional solution of
    \begin{equation}\label{eq:lem}
        -\left(\Delta + k^2 q(x)\right)v = g \ \text{in}\ \R^d.
    \end{equation}
    Then $v\in H^1$ and satisfies
    \begin{align}
        \|v\|_{L^2} &\le k^{-2}c^{-1}\|g\|_{L^2}, \label{eq:L2} \\
        \|\nabla v\|_{L^2}^2 &\le k^2\|q\|_{L^\infty}\|v\|_{L^2}^2 + \|g\|_{L^2}\|v\|_{L^2}, \label{eq:grad} \\
        \|v\|_{H^1} &\le k^{-1} C(c,\|q\|_{L^\infty}) \|g\|_{L^2}. \label{eq:H1}
    \end{align}
\end{lemma}
\begin{proof}
    Since the equation \eqref{eq:lem} implies $-\Delta v = k^2q v + g\in L^2$, we have $v\in H^2$. Testing the equation \eqref{eq:lem} with $v^*$ in the whole space, we obtain
    \begin{equation}\label{eq:id}
        \|\nabla v\|_{L^2}^2 - k^2\int_{\R^d} q|v|^2\,dx = \int_{\R^d} gv^*\,dx.
    \end{equation}
    Then \eqref{eq:grad} follows by direct estimate; taking the imaginary part gives
    \begin{equation*}
        k^2 c\|v\|_{L^2}^2 \le k^2\int_{\R^d} \Im q|v|^2\,dx = -\Im\int_{\R^d} gv^*\,dx \le \|g\|_{L^2} \|v\|_{L^2},
    \end{equation*}
    and \eqref{eq:L2} follows. Combining \eqref{eq:L2} and \eqref{eq:grad} implies \eqref{eq:H1}.
\end{proof}

One of the applications of the $L^2$ estimate \eqref{eq:L2} is to derive an estimate for the Green function of the homogenized equation \eqref{eq:homo}.
\begin{lemma}
    For any $y\in\R^d$, let $G_\epsilon(\cdot,y)$ be the solution of
    \begin{equation*}
        -\left(\Delta + k^2 \mu_\epsilon\right)G_\epsilon(\cdot,y) = \delta_y \ \text{in}\ \R^d
    \end{equation*}
    with $\mu_\epsilon$ defined by \eqref{eq:mu_eps}. Then for all $\epsilon\in(0,1)$,
    \begin{equation*}
        \sup_{y\in\R^d} \|G_\epsilon(\cdot,y)\|_{L^2}^2 \le C(k,\uns_0,\uns_1,d).
    \end{equation*}
\end{lemma}
\begin{proof}
    We firstly derive an estimate for the Green function $G_0$ of the operator $-\left(\Delta + k^2\uns_0\right)$ with constant coefficient. Since $\Im\uns_0 = \kappa_0 > 0$, we define $\uno$ to be the principal square root of $\uns_0$ (i.e.\ the evaluation of the single-valued branch of $\sqrt{\cdot}$ on $\mathbb{C}\setminus \R_+$), so that $\Im\uno > 0$. Then the Green function admits the following explicit formulas:
    \begin{equation*}
        G_0(x,y) = \begin{cases}
            \dfrac{e^{ik\uno|x-y|}}{4\pi|x-y|}, & d=3, \\[10pt]
            \dfrac{i}{4}H_0^{(1)}(k\uno|x-y|), & d=2,
        \end{cases}
    \end{equation*}
    where $H_0^{(1)}$ is the zeroth order Hankel function of the first kind.
    
    For $d=3$ and any $y\in\R^3$, by direct calculation we have
    \begin{align*}
        \|G_0(\cdot,y)\|_{L^2}^2 &= \int_{\R^3} |G_0(x,y)|^2\,dx \\
        &= \frac{1}{(4\pi)^2} \int_{\R^3} \frac{e^{-2k|x-y|\Im\uno}}{|x-y|^2}\,dx \\
        &= \frac{1}{4\pi} \int_0^\infty e^{-2kr\Im\uno}\,dr = \frac{1}{8\pi k\Im\uno}.
    \end{align*}
    For $d=2$, since the Hankel function $H_0^{(1)}$ has no analytic expressions, we need to involve its asymptotic behavior at infinity
    \begin{equation*}
        H_0^{(1)}(z) \sim \sqrt{\frac{2}{\pi z}} e^{i\left(z-\frac{\pi}{4}\right)}, \quad \arg z\in(-\pi,2\pi),\ |z|\to\infty.
    \end{equation*}
    Since
    \begin{align*}
        \int_{\R^2} \left|\sqrt{\frac{2}{\pi k\uno|x-y|}} e^{i\left(k\uno|x-y|-\frac{\pi}{4}\right)}\right|^2\,dx
        &= \int_{\R^2} \frac{2}{\pi k|\uno||x-y|} e^{-2k|x-y|\Im\uno}\,dx \\
        &= \frac{2}{\pi k|\uno|} \int_0^\infty e^{-2kr\Im\uno}\,dr \\
        &= \frac{1}{\pi k^2|\uno|\Im\uno},
    \end{align*}
    we obtain that
    \begin{equation*}
        \sup_{y\in\R^d} \|G_0(\cdot,y)\|_{L^2}^2 \le C(k,\uno).
    \end{equation*}

    Now we turn to prove the estimate of $G_\epsilon$. For any $y\in\R^d$, since
    \begin{equation*}
        -\left(\Delta + k^2\mu_\epsilon\right)(G_\epsilon(\cdot,y) - G_0(\cdot,y)) = k^2\left(\uns_1 - \uns_0\right)p_\epsilon G_0(\cdot,y) \ \text{in}\ \R^d
    \end{equation*}
    and
    \begin{equation*}
        \inf_{x\in\R^d} \Im\mu_\epsilon(x) = \inf_{x\in\R^d} \{p_\epsilon(x)\kappa_1 + (1-p_\epsilon(x))\kappa_0\} \ge \min\{\kappa_0,\kappa_1\} \eqqcolon \kappa_m > 0,
    \end{equation*}
    the $L^2$ estimate \eqref{eq:L2} gives
    \begin{align*}
        \|G_\epsilon(\cdot,y)-G_0(\cdot,y)\|_{L^2} &\le k^{-2}\kappa_m^{-1} \|k^2(\uns_1 - \uns_0)p_\epsilon G_0(\cdot,y)\|_{L^2} \\
        &\le \kappa_m^{-1} |\uns_1 - \uns_0| \|G_0(\cdot,y)\|_{L^2}.
    \end{align*}
    It follows that
    \begin{equation*}
        \|G_\epsilon(\cdot,y)\|_{L^2}^2 \le \left(1 + \kappa_m^{-1} |\uns_1 - \uns_0|\right)^2 \|G_0(\cdot,y)\|_{L^2}^2 \le C(k,\uns_0,\uns_1,d).
    \end{equation*}
\end{proof}

We are now ready to prove \cref{thm_homogenization}.
\begin{proof}
    First, we prove \eqref{eq:est1}. By \eqref{eq:hete} and \eqref{eq:inter} we have
    \begin{equation}\label{eq:diff}
        -(\Delta + k^2 \mu_\epsilon)(u_\epsilon - \tilde{u}_\epsilon) = k^2(\uns_\epsilon - \mu_\epsilon)u_\epsilon \ \text{in}\ \R^d.
    \end{equation}
    Denoting the solution operator of \eqref{eq:inter} by ${\cal G}_\epsilon$, which is the integral operator associated with the Green kernel $G_\epsilon(x,y)$, we have
    \begin{equation*}
        u_\epsilon - \tilde{u}_\epsilon = k^2{\cal G}_\epsilon r_\epsilon u_\epsilon
        = k^2{\cal G}_\epsilon r_\epsilon \tilde{u}_\epsilon + k^2{\cal G}_\epsilon r_\epsilon(u_\epsilon - \tilde{u}_\epsilon),
    \end{equation*}
    where $r_\epsilon(x) = \uns_\epsilon(x) - \mu_\epsilon(x)$ takes the value of $(1-p_\epsilon(x))(\uns_1-\uns_0)$ with probability $p_\epsilon(x)$ and $-p_\epsilon(x)(\uns_1-\uns_0)$ with probability $1-p_\epsilon(x)$.
    Define $v_\epsilon = k^2{\cal G}_\epsilon r_\epsilon \tilde{u}_\epsilon$, solving the equation
    \begin{equation}\label{eq:corrector}
        -(\Delta + k^2\mu_\epsilon)v_\epsilon = k^2 r_\epsilon \tilde{u}_\epsilon \ \text{in}\ \R^d.
    \end{equation}
    Then we need to estimate $v_\epsilon$ and $u_\epsilon - \tilde{u}_\epsilon - v_\epsilon$.
    For the latter, by \eqref{eq:hete}, \eqref{eq:inter} and \eqref{eq:corrector} we derive
    \begin{equation*}
        -(\Delta + k^2 \uns_\epsilon(x))(u_\epsilon - \tilde{u}_\epsilon - v_\epsilon) = k^2 r_\epsilon v_\epsilon \ \text{in}\ \R^d.
    \end{equation*}
    Since $\Im\uns_\epsilon(x)\ge \kappa_m$ for all $x\in\R^d$, the $L^2$ estimate \eqref{eq:L2} gives
    \begin{align*}
        \|u_\epsilon - \tilde{u}_\epsilon - v_\epsilon\|_{L^2} &\le k^{-2}\kappa_m^{-1} \|k^2 r_\epsilon v_\epsilon\|_{L^2} \\
        &\le \kappa_m^{-1} \|r_\epsilon\|_{L^\infty} \|v_\epsilon\|_{L^2}
        \le \kappa_m^{-1} |\uns_1 - \uns_0| \|v_\epsilon\|_{L^2}.
    \end{align*}
    So
    \begin{equation}\label{eq:diff_est}
        \|u_\epsilon - \tilde{u}_\epsilon\|_{L^2} \le (1 + \kappa_m^{-1} |\uns_1 - \uns_0|) \|v_\epsilon\|_{L^2}.
    \end{equation}
    Employing the Green function representation of $v_\epsilon$, we have
    \begin{align*}
        \|v_\epsilon\|_{L^2}^2 &= \int_{\R^d} v_\epsilon(x)v_\epsilon^*(x)\,dx \\
        &= k^4\int_{\R^d} \int_{\R^d} G_\epsilon(x,y)r_\epsilon(y)\tilde{u}_\epsilon(y)\,dy \int_{\R^d} G_\epsilon^*(x,z)r^*_\epsilon(z) \tilde{u}_\epsilon^*(z)\,dz\,dx,
    \end{align*}
    and it follows that
    \begin{align*}
        \be\|v_\epsilon\|_{L^2}^2 &= k^4 \int_{\R^d} \int_{\R^d\times\R^d} G_\epsilon(x,y) G_\epsilon^*(x,z) \be[r_\epsilon(y)r_\epsilon^*(z)] \tilde{u}_\epsilon(y) \tilde{u}_\epsilon^*(z)\,dy\,dz\,dx \\
        &= k^4 \int_{\R^d\times\R^d} \left(\int_{\R^d} G_\epsilon(x,y) G_\epsilon^*(x,z)\,dx\right) \be[r_\epsilon(y)r_\epsilon^*(z)] \tilde{u}_\epsilon(y) \tilde{u}_\epsilon^*(z)\,dy\,dz \\
        &\le C(k,\uns_0,\uns_1,d) \int_{\R^d\times\R^d} \big|\be[r_\epsilon(y)r_\epsilon^*(z)]\big| |\tilde{u}_\epsilon(y)||\tilde{u}_\epsilon(z)|\,dy\,dz.
    \end{align*}
    The right-hand integral reduces to the part where $|y-z|<2\epsilon$; it follows by the independence assumption of the Poisson point process and the definition of $r_\epsilon$. Indeed, when $|y-z|\ge 2\epsilon$, since $B(y,\epsilon)\cap B(z,\epsilon) = \varnothing$, the random variables $N_\epsilon(B(y,\epsilon))$ and $N_\epsilon(B(z,\epsilon))$ are independent. So are $r_\epsilon(y)$ and $r_\epsilon(z)$ being their function respectively. Hence,
    \begin{equation*}
        \be[r_\epsilon(y)r_\epsilon^*(z)] = \be[r_\epsilon(y)]\be[r_\epsilon(z)]^* = 0,
    \end{equation*}
    and the integral over $|y-z|\ge 2\epsilon$ vanishes.
    Then
    \begin{align*}
        \be\|v_\epsilon\|_{L^2}^2 \le{}& C(k,\uns_0,\uns_1,d)\int_{\{|y-z|<2\epsilon\}} \big|\be[r_\epsilon(y)r_\epsilon^*(z)]\big| |\tilde{u}_\epsilon(y)||\tilde{u}_\epsilon(z)|\,dy\,dz \\
        \le{}& C(k,\uns_0,\uns_1,d) \int_{\{|y-z|<2\epsilon\}} |\tilde{u}_\epsilon(y)||\tilde{u}_\epsilon(z)|\,dy\,dz \\
        ={}& C(k,\uns_0,\uns_1,d) \int_{\R^d} |\tilde{u}_\epsilon(y)| \int_{B(y,2\epsilon)} |\tilde{u}_\epsilon(z)|\,dz\,dy.
    \end{align*}
    Denoting
    \begin{equation*}
        {\cal M}\tilde{u}_\epsilon(y)  = \sup_{r>0} \fint_{B(y,r)} |\tilde{u}_\epsilon(z)|\,dz
    \end{equation*}
    as the Hardy--Littlewood maximal function of $\tilde{u}_\epsilon$, we further have
    \begin{align*}
        \be\|v_\epsilon\|_{L^2}^2 &\le C(k,\uns_0,\uns_1,d) V_d(2\epsilon)^d \int_{\R^d} |\tilde{u}_\epsilon(y)||{\cal M}\tilde{u}_\epsilon(y)|\,dy \\
        &\le C(k,\uns_0,\uns_1,d) \epsilon^d \|\tilde{u}_\epsilon\|_{L^2} \|{\cal M}\tilde{u}_\epsilon\|_{L^2}.
    \end{align*}
    Incorporated with the maximal function theorem \cite[Theorem 2.2]{Heinonen2001}
    \begin{equation*}
        \|{\cal M}\tilde{u}_\epsilon\|_{L^2} \le C(d)\|\tilde{u}_\epsilon\|_{L^2}
    \end{equation*}
    and the $L^2$ estimate \eqref{eq:L2} for $\tilde{u}_\epsilon$, it follows that
    \begin{equation*}
        \be\|v_\epsilon\|_{L^2}^2 \le C(k,\uns_0,\uns_1,d) \epsilon^d \|\tilde{u}_\epsilon\|_{L^2}^2 \le C(k,\uns_0,\uns_1,d) \epsilon^d \|f\|_{L^2}^2.
    \end{equation*}
    Together with \eqref{eq:diff_est}, this proves \eqref{eq:est1}.

    Now we continue to prove \eqref{eq:est2}. By \eqref{eq:inter} and \eqref{eq:homo} we have
    \begin{equation}\label{eq:diff_deter}
        -(\Delta + k^2\mu_\epsilon)(\tilde{u}_\epsilon - u) = k^2(\mu_\epsilon - \mu)u.
    \end{equation}
    Applying the $L^2$ estimate \eqref{eq:L2}, we derive
    \begin{equation}\label{eq:diff_est2}
        \|\tilde{u}_\epsilon - u\|_{L^2} \le \kappa_m^{-1}\|\mu_\epsilon - \mu\|_{L^\infty} \|u\|_{L^2} \le C(k,\uns_0,\uns_1)\|p_\epsilon - p\|_{L^\infty} \|f\|_{L^2}.
    \end{equation}
    Since $\rho$ is nonnegative and $C^\alpha$, we have
    \begin{align*}
        |p_\epsilon(x) - p(x)| &= \big|e^{-\fint_{B(x,\epsilon)}\rho(y)\,dy} - e^{-\rho(x)}\big| \\
        &\le \bigg|\fint_{B(x,\epsilon)} \rho(y)\,dy - \rho(x)\bigg| \\
        &\le \fint_{B(x,\epsilon)} |\rho(y) - \rho(x)|\,dy 
        \le C(d,\alpha,[\rho]_\alpha) \epsilon^\alpha
    \end{align*}
    for all $x\in\R^d$. Along with \eqref{eq:diff_est2}, this proves \eqref{eq:est2}.
    
    Finally, since $2\alpha\le 2\le d$, \eqref{eq:homo_est} follows from \eqref{eq:est1} and \eqref{eq:est2}.
\end{proof}

We conclude this subsection with the proof of \cref{cor:homo}.
\begin{proof}
    Applying the gradient estimate \eqref{eq:grad} to the equation \eqref{eq:diff} yields that
    \begin{align*}
        \|\nabla(u_\epsilon - \tilde{u}_\epsilon)\|_{L^2}^2 
        &\le k^2\|\mu_\epsilon\|_{L^\infty} \|u_\epsilon - \tilde{u}_\epsilon\|_{L^2}^2 + k^2\|r_\epsilon\|_{L^\infty} \|u_\epsilon\|_{L^2} \|u_\epsilon - \tilde{u}_\epsilon\|_{L^2} \\
        &\le C(k,\uns_0,\uns_1) \left(\|u_\epsilon - \tilde{u}_\epsilon\|_{L^2}^2 + \|f\|_{L^2} \|u_\epsilon - \tilde{u}_\epsilon\|_{L^2}\right).
    \end{align*}
    Taking the expectation, by the Cauchy--Schwarz inequality and the estimate \eqref{eq:est1}, we derive
    \begin{align*}
        \be\|\nabla(u_\epsilon - \tilde{u}_\epsilon)\|_{L^2}^2
        &\le C(k,\uns_0,\uns_1) \left(\be\|u_\epsilon - \tilde{u}_\epsilon\|_{L^2}^2 + \|f\|_{L^2} \be\|u_\epsilon - \tilde{u}_\epsilon\|_{L^2}\right) \\
        &\le C(k,\uns_0,\uns_1) \left(\be\|u_\epsilon - \tilde{u}_\epsilon\|_{L^2}^2 + \|f\|_{L^2} \big(\be\|u_\epsilon - \tilde{u}_\epsilon\|_{L^2}^2\big)^{1/2}\right) \\
        &\le C(k,\uns_0,\uns_1,d) \epsilon^{d/2} \|f\|_{L^2}^2.
    \end{align*}
    Then the estimate \eqref{eq:est1_H1} follows.
    
    For the estimate \eqref{eq:est2_H1}, we only have to apply the $H^1$ estimate \eqref{eq:H1} to the equation \eqref{eq:diff_deter}. The boundary estimate \eqref{eq:homo_est_bdry} then follows by combining \eqref{eq:est1_H1} and \eqref{eq:est2_H1} and employing the trace theorem.
\end{proof}

\begin{remark}
    There is an order loss in the random fluctuation error when passing the $L^2$ estimate to the $H^1$ estimate. It is caused by the energy estimate for the gradient, which contains one factor of the $L^2$ error. The deterministic bias retains the rate $\epsilon^\alpha$ in $H^1$, since its forcing term is already of order $\epsilon^\alpha$ in $L^2$.
\end{remark}

\subsection{Multiscale inverse problem}
In the following, we denote the data generated by the heterogeneous medium as $d_{\rm het}$ and the data generated by the effective medium as $d_{\rm eff}$. Due to the expensive computational cost of the direct inversion of the operator ${\cal F}_{\rm het}^\epsilon$, many frameworks have been proposed to solve the inverse problems associated with the multiscale nature. Inspired by the study in \cite{frederickNumericalMethodsMultiscale2017} and considering the explicit form of the homogenized equation, we propose the following multiscale inverse problem formulation:
\begin{equation}\label{eq_inverse_problem}
    \min_{m(x)}\dist({\cal F}_{\rm eff}(m),d_{\rm het}),
\end{equation}
which replaces the synthetic data $d_{\rm eff}$ with the real measurements $d_{\rm het}$. Here the distance function $\dist(\cdot,\cdot)$ is a suitable metric to measure the discrepancy between the synthetic data and the real measurements.

The formulation \eqref{eq_inverse_problem} bears a strong resemblance to the full waveform inversion (FWI, \cite{laillySeismicInverseProblem1983,tarantolaInversionSeismicReflection1984a}) problem in seismology, which seeks to reconstruct the subsurface model using seismic data. However, there is a fundamental difference between the two. In FWI, the observed data is typically assumed to lie within the range of the forward operator, disregarding the effects of noise or model imperfections. In contrast, the data $d_{\rm het}$ in our problem exhibits multiscale and high-frequency characteristics, introducing complexities that prevent it from being contained within the range of ${\cal F}_{\rm eff}$. To address these complexities, we summarize the key challenges of solving this multiscale inverse problem as follows:
\begin{enumerate}
    \item The operator ${\cal F}_{\rm eff}$ is nonlinear, which makes the inverse problem ill-posed and difficult to solve numerically.
    \item The data $d_{\rm het}$ has a multiscale, stochastic and high-frequency nature, which makes the inverse problem underdetermined.
\end{enumerate}

To address these challenges, we propose the following numerical strategies:
\begin{enumerate}
    \item A coarse-grid discretization of the effective medium $m(x)$ to reduce the dimensionality of the optimization problem and improve computational efficiency.
    \item Combination of measurements corresponding to different
    realizations of the microscopic medium through shot averaging, in order to diminish the impact of the random fluctuation of the data and enhance the robustness of the inversion.
    \item Mollification of the effective medium during the inversion to suppress oscillatory artifacts and incorporate the expected macroscopic smoothness of the reconstructed coefficient.
\end{enumerate}
These strategies and their numerical implementation are discussed in the following section.

\section{Numerical approach}\label{sec_numerical}
\begin{figure}[t]
    \centering
    \begin{tikzpicture}[auto,node1/.style={rectangle,draw=black, thick,fill=gray!20,text width=4em, text centered, rounded corners},arrow1/.style={thick,->,>=stealth}]
        \matrix[column sep=15mm, row sep=10mm]{
            \node[node1] (p) { $p(x)$}; & \node[node1] (c) {$c_\epsilon(x)$ }; & \node[node1] (rec) { $d_{\rm het}$}; \\
            & \node[node1] (homo_c){$c(x)$}; & \\
            \node[node1 ] (recoveryp) { $\tilde{p}(x)$}; & \node[node1 ] (recoveryc) {$\tilde{c}(x)$}; & \node[node1] (rec_coarse) { $d_{\rm het}$}; \\
        };
        \node[below=of recoveryc,text=gray](null){Direct reconstruction};

        \draw[arrow1] (p)--(c);
        \draw[arrow1] (c)--(rec);
        \draw[arrow1] (p)|-(homo_c);
        \draw[arrow1] (c)--node[right,near end]{Homogenization}(homo_c);
        \draw[arrow1] (rec)--node[right]{Down-sampling along $t$ axis}(rec_coarse);
        \draw[arrow1] (rec_coarse)--node[above]{FWI}(recoveryc);
        \draw[arrow1] (recoveryc)--(recoveryp);
        \draw[arrow1,gray,dashed] (rec_coarse)|-(null.south)-|(recoveryp);
        \draw[thick,dashed] (recoveryc)--node[right,near end]{Approximation}(homo_c);

        \node[draw=red,text=red, dotted,thick,inner sep=1em,fit=(c) (rec),label={[text=red]above:Forward model (fine mesh)}]{};
        \node[draw=blue,text=blue, dotted,thick,inner sep=1em,fit=(rec_coarse) (recoveryc) (recoveryp) (null),label={[text=blue]below:Inverse problem (coarse mesh)}]{};
    \end{tikzpicture}
    \vspace*{8pt}
    \caption{The diagram of numerical simulations.}
    \label{fig_diagram}
\end{figure}
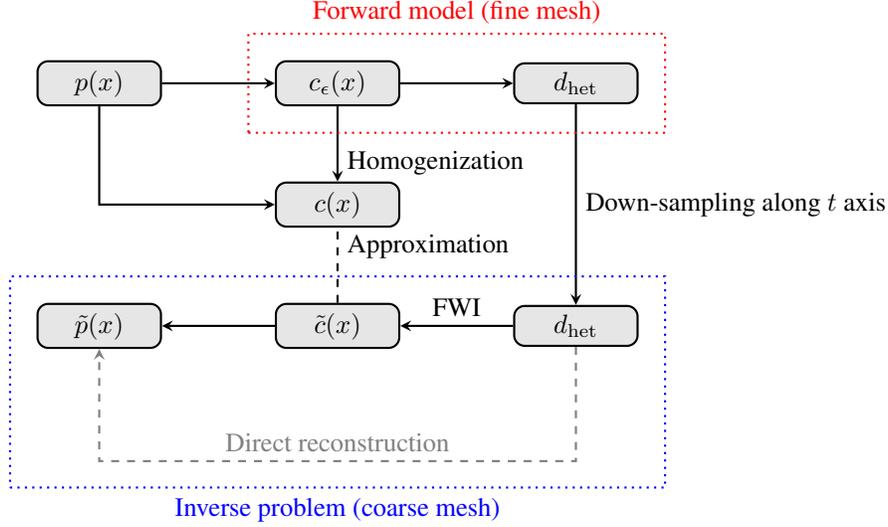

The overall workflow of the numerical simulations is illustrated in \cref{fig_diagram}. Broadly, the process consists of two components: forward modeling (fine mesh) and solving the inverse problem (coarse mesh). The forward modeling involves generating synthetic data $d_{\rm het}$ based on the heterogeneous medium $c_\epsilon(x)$, while the inverse problem focuses on recovering the effective medium $c(x)$ and the sediment distribution $p(x)$ from the measurements $d_{\rm het}$. It is important to note that the synthetic data $d_{\rm het}$ in forward modeling is generated on a fine mesh to mimic real measurements and capture the multiscale characteristics of the data. In contrast, all simulations for the inverse problem are carried out on a coarse mesh. Once we get the heterogeneous data $d_{\rm het}$, we will first preprocess these data to obtain the data suitable for inversion, and then solve the inverse problem to get $\tilde{c}(x)$ and $\tilde{p}(x)$. According to \cref{thm_homogenization}, the recovered effective medium $\tilde{c}(x)$ is expected to closely match the true effective medium $c(x)$ corresponding to the heterogeneous medium $c_\epsilon(x)$. Besides, one could directly reconstruct the sediment distribution $\tilde{p}(x)$ from the data $d_{\rm het}$ skipping the medium recovery step. In the following, we will introduce the numerical methods for each step of the diagram in detail. 

\subsection{Effective model and its parameterization}

For the physical modeling of the sediment distribution in \cref{sec_physical_model_sediment}, we have derived the effective medium model $m(x)$ according to \eqref{eq_m_eff}, which is
\begin{equation}
    m(x) = \frac{1}{c^2(x)} = \frac{p(x)}{c_1^2} + \frac{1-p(x)}{c_0^2}.
\end{equation}

For the numerical inversion, the effective medium is discretized on the coarse finite-difference grid. The optimization problem \eqref{eq_inverse_problem} is accordingly discretized and the grid values of the effective medium are taken as the optimization variables. This coarse-grid parameterization is consistent with our objective of recovering only the macroscopic effective medium rather than the microscopic variations of the heterogeneous coefficient.

\begin{remark}
    In the present numerical experiments, all values of the effective medium on the coarse finite-difference grid are treated as independent optimization variables. A further reduction of the parameter dimension could be achieved by representing the effective medium in a suitable reduced basis, or by employing a sparse representation when additional prior information on the macroscopic sediment profile is available. Such reduced parameterizations are left for future work.
\end{remark}

\subsection{Methods to solve the inverse problems}\label{sec_solveIP}
To solve the inverse problem \eqref{eq_inverse_problem}, gradient-based optimization methods are employed to minimize the misfit between the observed and predicted data. The predicted data is generated by solving the wave equation, and the gradient of the misfit is calculated by solving the adjoint wave equation. The choice of objective function plays a critical role in influencing the quality of inversion results. In this study, we consider two objective functions: the $L^2$ norm and the quadratic Wasserstein distance $W^2$. Detailed discussions on FWI based on these objective functions can be found in \cite{engquistOptimalTransportBased2022,plessixReviewAdjointstateMethod2006,qiuAnalysisSeismicInversion2021a}. The $W^2$ objective function is particularly appealing for several reasons:
\begin{enumerate}
    \item  It facilitates convergence even from a rough initial guess, which is crucial for inverse problems with limited prior knowledge.
    \item  It promotes smoothness in the recovered model, which aligns well with the smooth nature of the target effective model in our inversion problem.
    \item  It is demonstrated that $W^2$ is not very sensitive to oscillations, as highlighted in \cite{qiuAnalysisSeismicInversion2021a}, making it well-suited for applications involving multiscale features.
\end{enumerate}

Next, we introduce a shot averaging strategy, or more precisely a realization averaging strategy, to reduce the influence of the microscopic random configuration on the inversion. The estimates in \cref{thm_homogenization} and \cref{cor:homo} distinguish the realization-dependent fluctuation error from the deterministic scale bias. The purpose of the averaging procedure is to suppress the former contribution before solving the inverse problem. More precisely, we generate $N$ independent realizations $c_\epsilon^{(1)}(x),\ldots,c_\epsilon^{(N)}(x)$ from the same macroscopic probability pattern $p(x)$, and solve the heterogeneous wave equation for each realization to obtain the corresponding boundary measurements $d_{\rm het}^{(1)},\ldots,d_{\rm het}^{(N)}$. The averaged data are defined by
\begin{equation}
    \bar{d}_{\rm het} \coloneqq \frac{1}{N}\sum_{j=1}^N d_{\rm het}^{(j)}.
    \label{eq_shot_averaging}
\end{equation}
Then we replace the synthetic data $d_{\rm het}$ with the averaged data $\bar{d}_{\rm het}$ in the inverse problem \eqref{eq_inverse_problem}. The independence assumption of the realizations corresponds physically to repeated measurements taken from microscopic sediment configurations that have sufficiently decorrelated between successive acquisitions, while the macroscopic probability pattern $p(x)$ remains unchanged.

The fluctuation-reduction effect of this procedure can be quantified in the Hilbert data space ${\cal H} \coloneqq L^2(\Gamma\times(0,T))$. For independent realizations, we have
\begin{equation*}
    \be\big\|\bar{d}_{\rm het} - \be[d_{\rm het}]\big\|_{\cal H}^2 = \frac{1}{N}\be\big\|d_{\rm het} - \be[d_{\rm het}]\big\|_{\cal H}^2.
\end{equation*}
Thus, averaging reduces the variance of the realization-dependent part by a factor of $N^{-1}$. However, it does not eliminate the systematic difference between the mean heterogeneous data and the effective-model data at a fixed scale $\epsilon>0$.

For the $L^2$ objective, averaging the data before inversion has an additional exact interpretation. Indeed, for any candidate effective medium $m$,
\begin{equation*}
    \frac{1}{N}\sum_{j=1}^N \big\|{\cal F}_{\rm eff}(m)-d_{\rm het}^{(j)}\big\|_{\cal H}^2
    = \big\|{\cal F}_{\rm eff}(m)-\bar{d}_{\rm het}\big\|_{\cal H}^2 
    + \frac{1}{N}\sum_{j=1}^N \big\|d_{\rm het}^{(j)}-\bar{d}_{\rm het}\big\|_{\cal H}^2.
\end{equation*}
Since the second term is independent of $m$, minimizing the $L^2$ misfit against the averaged data is equivalent to minimizing the average of the individual $L^2$ misfits.

For the $W^2$ objective, however, no analogous identity generally holds. In that case, replacing the individual measurements by their average $\bar{d}_{\rm het}$ should be regarded as a preprocessing choice: it is intended to attenuate the realization-dependent fluctuations and to produce data that better approximate the macroscopic wave response before the inversion is performed. We do not claim that minimizing the $W^2$ misfit against the averaged data is equivalent to averaging the individual $W^2$ misfits; these are different optimization formulations. The numerical experiments in \cref{sec_results} assess the performance of the averaged-data formulation for both the $L^2$ and $W^2$ objectives.

\begin{remark}
    If the probability distribution $p(x)$ exhibits a low-rank parameterization, one might contemplate bypassing FWI altogether and directly inverting to obtain $p(x)$. This alternative could potentially offer a more stable and efficient approach, as it involves a reduction in the dimensionality of the unknown parameter space.
\end{remark}

The shot averaging strategy reduces the fluctuations associated with different microscopic realizations at the data level. A complementary strategy is to exploit the macroscopic smoothness of the effective medium itself. For this purpose, we consider a model mollification strategy: let $K$ be a mollifying kernel and consider
\begin{equation}
    \min_{m(x)}\dist({\cal F}_{\rm eff}(K\ast m),d_{\rm het}).
\end{equation}
The convolution $K\ast m$ suppresses cell-scale oscillations in the optimization variables and promotes a smoother effective coefficient, which is consistent with the macroscopic nature of the target medium.
In the numerical implementation, the mollification is carried out through the corresponding discrete convolution on the finite-difference grid.

Thus, shot averaging and model mollification act on two different aspects of the multiscale inverse problem: the former reduces realization-dependent fluctuations in the measurements, whereas the latter suppresses oscillatory artifacts in the reconstructed macroscopic coefficient. They can therefore be used either separately or in combination.

\begin{remark}
    Other stabilization strategies can also be incorporated into the present framework. For example, one may mollify both the predicted and observed data:
    \begin{equation}
        \min_{m(x)}\dist(K\ast({\cal F}_{\rm eff}(m)),K\ast d_{\rm het});
    \end{equation}
    or introduce an explicit regularization term,
    \begin{equation}
        \min_{m(x)}\dist({\cal F}_{\rm eff}(m),d_{\rm het})+\alpha \mathcal{R}(m_H).
    \end{equation}
    The regularization term $\mathcal{R}$ is introduced to prevent overfitting and stabilize the solution by imposing additional constraints on the model parameters. It may be chosen, for example, as a Tikhonov, total variation or sparsity-promoting regularization term. These alternatives are not investigated in the numerical experiments of the present work.
\end{remark}

\subsection{Determination of sediment concentration}
Once we obtain the effective medium model $\tilde{c}(x)$, we can determine the sediment distribution $\tilde{p}(x)$ through direct calculation: 
\begin{equation} 
    \tilde{p}(x)=\frac{\frac{1}{\tilde{c}^2(x)}-\frac{1}{c_0^2}}{\frac{1}{c_1^2}-\frac{1}{c_0^2}}. 
\end{equation} 
Subsequently, the estimated sediment concentration can be determined by integrating $\tilde{p}(x)$ over the domain $D$. Note that the sediment concentration calculated here is the ratio of the volume of the sediment to the volume of the water, rather than the conventional mass concentration, i.e.\ the ratio of the mass of the sediment to the volume of the water. However, the mass concentration can be easily calculated by multiplying the volume concentration by the density of the sediment. In the following numerical experiments, we will focus on the volume concentration.

\section{Numerical results}\label{sec_results}
In this section, we present the numerical results to demonstrate the effectiveness of the proposed model and methods. 

The numerical experiments are conducted in a two-dimensional domain $[0,\SI{1}{\meter}]\times[0,\SI{1}{\meter}]$. 
The heterogeneous medium is generated on a $5000\times5000$ fine grid with spacing
\begin{equation*}
    h = \SI{2e-4}{\meter},
\end{equation*}
where each fine cell is assigned either the sediment or the water coefficient according to the prescribed probability pattern. The inverse problem is solved on a $500\times500$ coarse grid with spacing
\begin{equation*}
    H=\SI{2e-3}{\meter}=10h,
\end{equation*}
and the values of the effective coefficient at the coarse-grid points are taken as the optimization variables.

\begin{remark}
    The cellwise binary construction is used as a discrete surrogate of the Poisson-cloud medium at the microscopic scale. It preserves the binary heterogeneous structure and the prescribed local sediment-occurrence probability, although it does not reproduce the spherical particle geometry exactly. For the purpose of scale comparison, interpreting the fine-cell width as the particle diameter gives $\epsilon\approx h/2$, and hence $\epsilon/H\approx 0.05$.
\end{remark}

The p-wave velocity of water is \(c_0=\SI[per-mode=symbol]{1500}{\meter\per\second}\), and that of sand is \(c_1=\SI[per-mode=symbol]{3000}{\meter\per\second}\). The source (red dots in the following figures) is located along the bottom side of the domain, while the receivers (green diamonds) are deployed along the left, top, and right sides. The source is a Ricker wavelet with a central frequency of \(\SI{15}{\kilo\hertz}\), which is treated as effectively band-limited over the frequency range resolved by the numerical discretization, since its spectral energy outside this range is negligible. The time steps are \(\SI{20}{\nano\second}\) and \(\SI{500}{\nano\second}\) for the fine and coarse meshes, respectively. The total simulation time is \(\SI{1}{\milli\second}\).

\begin{figure}[t]
    \centerline{\includegraphics[width=0.8\textwidth]{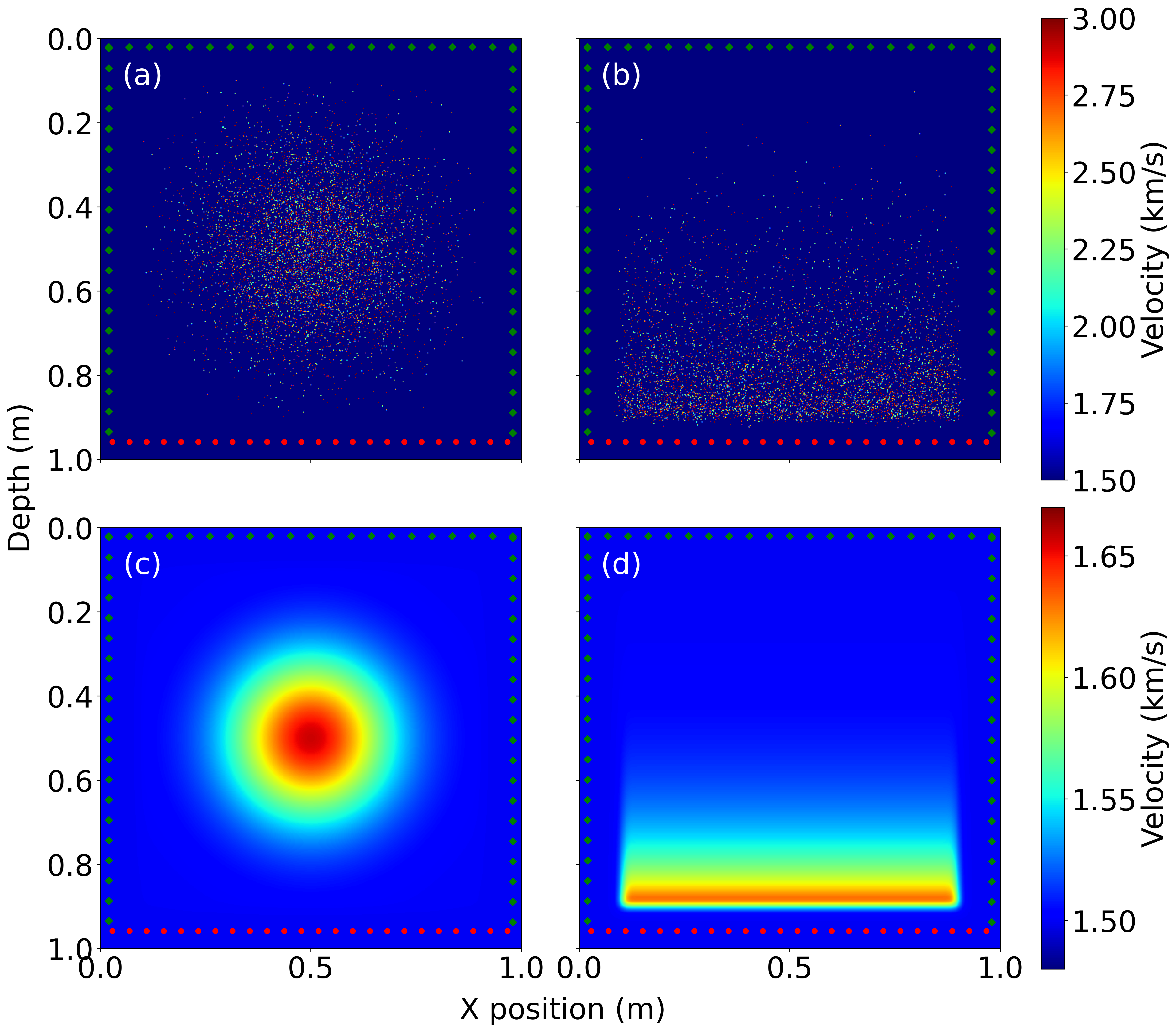}}
    \vspace*{8pt}
    \caption{The sediment distribution. (a) Sediment sample $c_\epsilon(x)$ from Gaussian function. (b) Sediment sample $c_\epsilon(x)$ from \cite{chiuMathematicalModelsDistribution2000}. (c) Effective medium model $c(x)$ from Gaussian function. (d) Effective medium model $c(x)$ from \cite{chiuMathematicalModelsDistribution2000}.}
    \label{fig_sediment_concentration}
\end{figure}

For the macroscopic probability field $p(x)$, we consider two cases: a Gaussian function and an empirical profile inferred from the real data \cite{chiuMathematicalModelsDistribution2000}.
The Gaussian function is given by
\begin{equation}
    p(x)=p_{\max}\exp\left(-\frac{|x-x_0|^2}{2\sigma^2}\right),
    \label{eq_gaussian}
\end{equation}
and the empirical profile inferred from the real data is given by 
\begin{equation}
    p(x)=p_{\max}\left(\frac{x_2}{e^M-x_2(e^M-1)}\right)^{\lambda G},~G=\frac{1-e^{-M}}{M\phi},~\phi=\frac{e^M}{e^M-1}-\frac{1}{M},
    \label{eq_chiu}
\end{equation}
where $p_{\max}$ is the peak probability, $x=(x_1,x_2)$ is the spatial variable, $x_2$ represents the depth, $M$ and $\lambda$ are the parameters. The parameters in the following experiments are chosen as $p_{\max}=0.25$, $x_0=(\SI{0.5}{\meter},\SI{0.5}{\meter})$, $\sigma=700$, $M=1$, and $\lambda=2$. The sediment distribution is shown in \cref{fig_sediment_concentration}. The optimization is performed by the L-BFGS-B algorithm.

\subsection{Verification of the shot data containing information of sediment concentration}
In this section, we first verify that the shot data contain information about the sediment concentration. We compare the shot data in the sediment-laden and sediment-free media. The shot data in the sediment-free media is generated by the homogeneous medium with the p-wave velocity $c_0=\SI[per-mode=symbol]{1500}{\meter\per\second}$. The shot data in the sediment-laden media is generated by the heterogeneous medium with the macroscopic probability field $p(x)$ given by \eqref{eq_gaussian} and \eqref{eq_chiu}. 

\begin{figure}[t]
    \centerline{\includegraphics[width=0.9\textwidth]{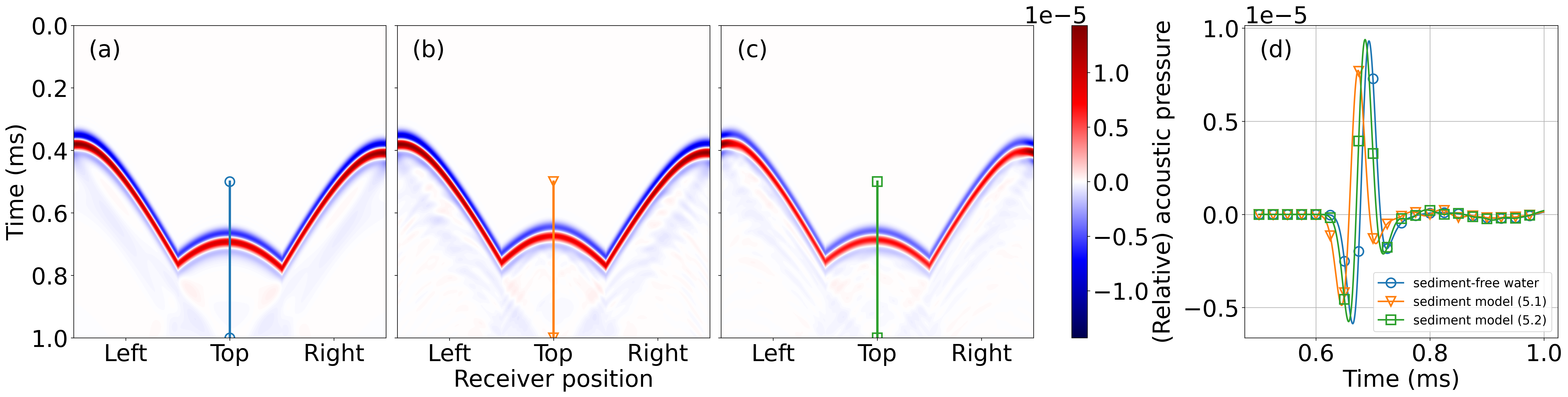}}
    \vspace*{8pt}
    \caption{Comparison of shot data in the sediment-laden and sediment-free media. (a) Shot data in the sediment-free media. (b) Shot data in the sediment-laden model \eqref{eq_gaussian}. (c) Shot data in the sediment-laden model \eqref{eq_chiu}. (d) Comparison of the shot data in the sediment-laden and sediment-free media.}
    \label{fig_comparison}
\end{figure}

The measurements are shown in \cref{fig_comparison}. \Cref{fig_comparison}(a) shows the shot data in the sediment-free media, while \cref{fig_comparison}(b) and \cref{fig_comparison}(c) show the shot data in the sediment-laden media with the macroscopic probability field given by \eqref{eq_gaussian} and \eqref{eq_chiu}, respectively. \Cref{fig_comparison}(d) presents the comparison of the shot data in the sediment-laden and sediment-free media.

It can be observed that the shot data in the sediment-laden media exhibit a different pattern compared to the shot data in the sediment-free media, especially the first arrival. The difference in the shot data is due to the presence of sediment, which affects the wave propagation in the medium. This observation confirms that the shot data contain information about the sediment concentration, which can be utilized to estimate the sediment concentration through the inverse problem.

\subsection{Data comparison}
To demonstrate the approximation of the effective medium model, we first compare the measurements of the wave field in the heterogeneous medium and the effective medium. In this analysis, data $d_{\rm het}$ is generated using a heterogeneous medium on a fine grid, while simultaneously $d_{\rm eff}$ is produced using the associated effective medium on a coarser grid. 
The receivers are deployed over the same physical locations (evenly distributed over the left, top and right sides). 

\begin{figure}[t]
    \centerline{\includegraphics[width=0.9\textwidth]{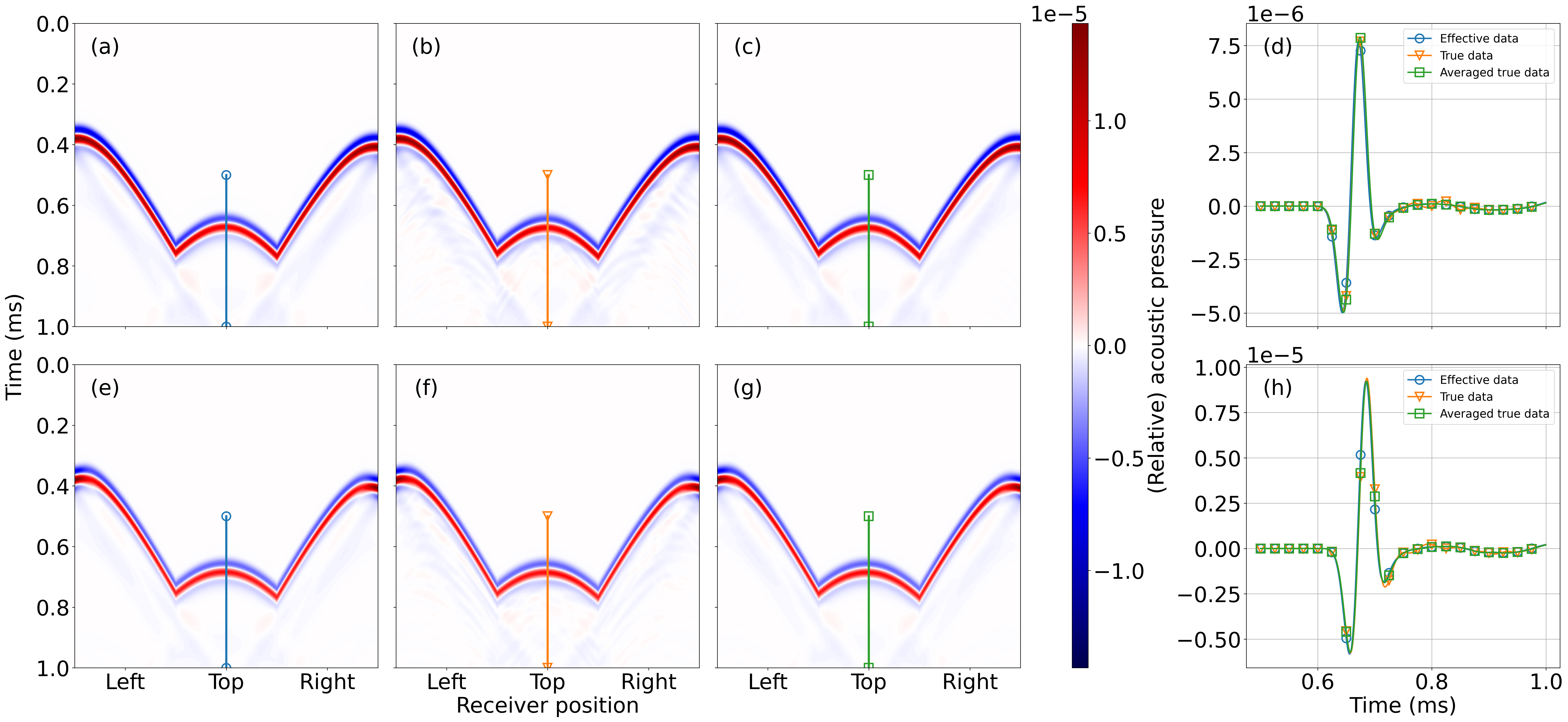}}
    \vspace*{8pt}
    \caption{The measurements of the wave field in the heterogeneous medium and the effective medium. The top row shows the results for the Gaussian function, while the bottom row shows the results for the sediment profile from \cite{chiuMathematicalModelsDistribution2000}. (a) Effective data $d_{\rm eff}$. (b) Heterogeneous data $d_{\rm het}$. (c) Averaged heterogeneous data. (d) Comparison for the data from the Gaussian function. (e) Effective data $d_{\rm eff}$. (f) Heterogeneous data $d_{\rm het}$. (g) Averaged heterogeneous data. (h) Comparison for the data from \cite{chiuMathematicalModelsDistribution2000}.}
    \label{fig_measurement}
\end{figure}

The measurements are shown in \cref{fig_measurement}. The top row in \cref{fig_measurement} shows the results for the Gaussian function, while the bottom row shows the results for the sediment profile from \cite{chiuMathematicalModelsDistribution2000}. The first column shows the effective data $d_{\rm eff}$, the second column shows the heterogeneous data $d_{\rm het}$, and the third column shows the averaged heterogeneous data. The last column shows the comparison of the data in the corresponding colored lines.

It can be observed that the measurements of the wave field in the effective medium are a good approximation of the measurements in the heterogeneous medium except for the multiple scattering parts. In particular, there is only a minor difference between the effective data and the averaged heterogeneous data. The mismatch of the multiple reflection parts is reasonable because the wave field generated by the effective medium model is a leading-order approximation of the wave field generated by the heterogeneous medium and the minor oscillations align with the high-order terms in the homogenized expansions. For a more accurate approximation, one can consider the higher-order homogenization theory, which can be found in \cite{capdevilleSecondOrderHomogenization2007,chenDispersiveModelWave2000,fishHigherOrderHomogenizationInitial2001,fishNonlocalDispersiveModel2002a}.
Overall, the effective medium model captures the direct wave quite well, which is the most important part of the wave field. This result also provides numerical validation of the effectiveness of the effective medium model claimed in \cref{thm_homogenization}.

\subsection{Inverse crime}
Utilizing the effective medium, we generate the synthetic data and perform the inversion with it by FWI. Equivalently, we solve the following minimization problem:
\begin{equation}
    \min_{m(x)}\dist\left({\cal F}_{\rm eff}(m),d_{\rm eff}\right).
\end{equation}
For the inverse crime, all the computations are performed on the coarse mesh, which is optimal for the effective medium. 

\begin{figure}[t]
    \centerline{\includegraphics[width=0.8\textwidth]{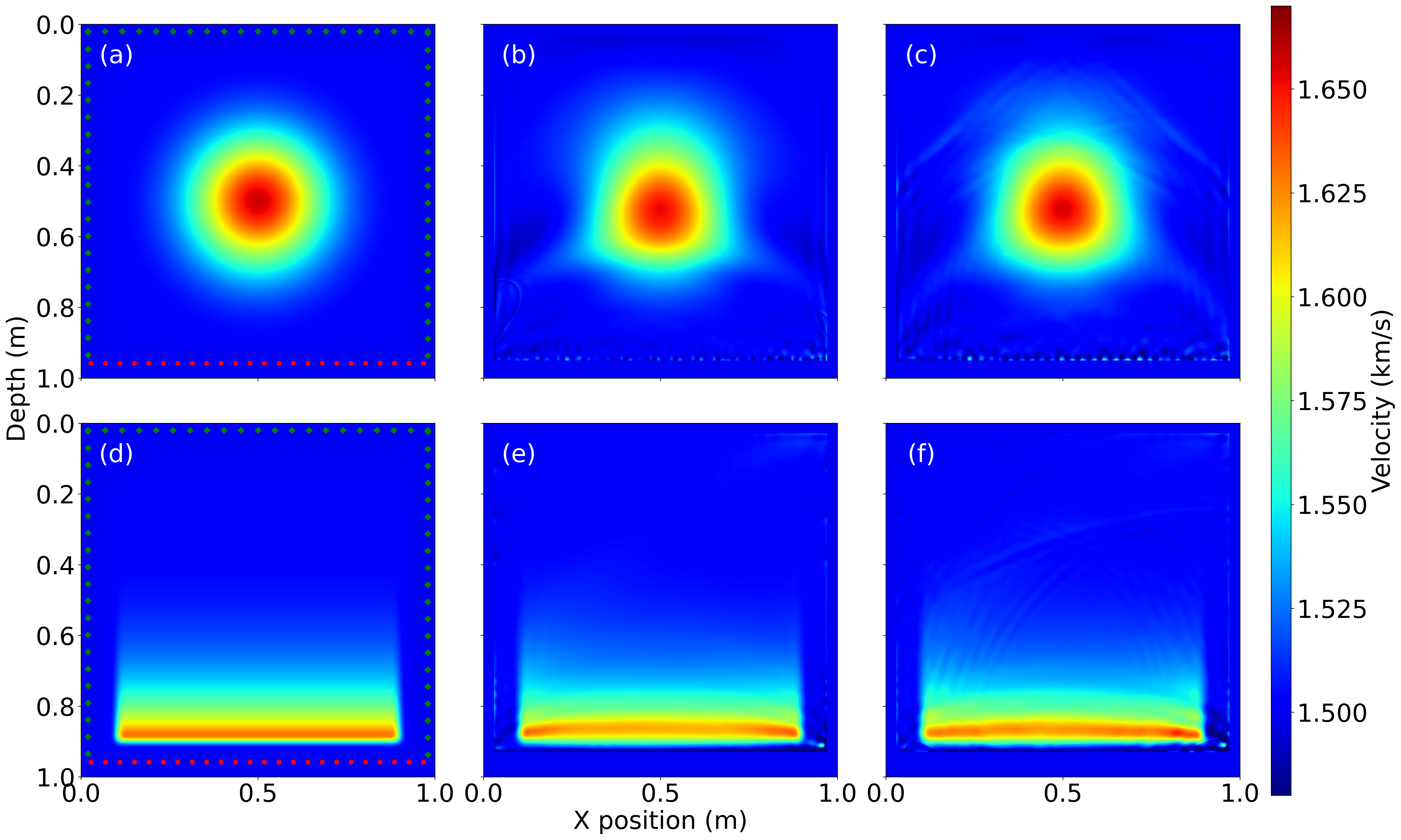}}
    \vspace*{8pt}
    \caption{The true model, and inverted models with $L^2$ and $W^2$ objective function. The top row shows the results for the Gaussian function, while the bottom row shows the results for the sediment profile from \cite{chiuMathematicalModelsDistribution2000}. (a) True effective model. (b) Inverted model with $L^2$ objective function. (c) The inverted model with $W^2$ objective function. (d) True effective model. (e) Inverted model with $L^2$ objective function. (f) The inverted model with $W^2$ objective function.}
    \label{fig_inverse_crime}
\end{figure}

In the following, we present the true model and inverted models with $L^2$ and $W^2$ objective functions in \cref{fig_inverse_crime}. The top row in \cref{fig_inverse_crime} shows the results for the Gaussian function, while the bottom row shows the results for the sediment profile from \cite{chiuMathematicalModelsDistribution2000}. The left column shows the true effective model, the middle column shows the inverted model with the $L^2$ objective function, and the right column shows the inverted model with the $W^2$ objective function.

It is noteworthy that both the $L^2$ and $W^2$ approaches produce effective reconstructions in this `inverse crime' scenario. The model using the $W^2$ objective function shows a slightly smaller error. However, the $L^2$ inverted model yields a smoother result, while the $W^2$ inverted model exhibits some unexpected oscillatory artifacts near the boundary. These artifacts can be effectively mitigated through the model mollification strategy introduced in \cref{sec_solveIP}.

\subsection{Inversion with data generated with heterogeneous medium}
In this numerical experiment, we employ the heterogeneous medium to generate the synthetic data and subsequently perform the inversion with it. More precisely, we solve the minimization problem \eqref{eq_inverse_problem}. For the shot averaging technique, we generate $N=50$ realizations of the heterogeneous medium using the same probability pattern $p(x)$ and solve the wave equation for each realization to obtain the corresponding measurements. The averaged data is then computed according to \eqref{eq_shot_averaging}.

\begin{figure}[t]
    \centerline{\includegraphics[width=0.8\textwidth]{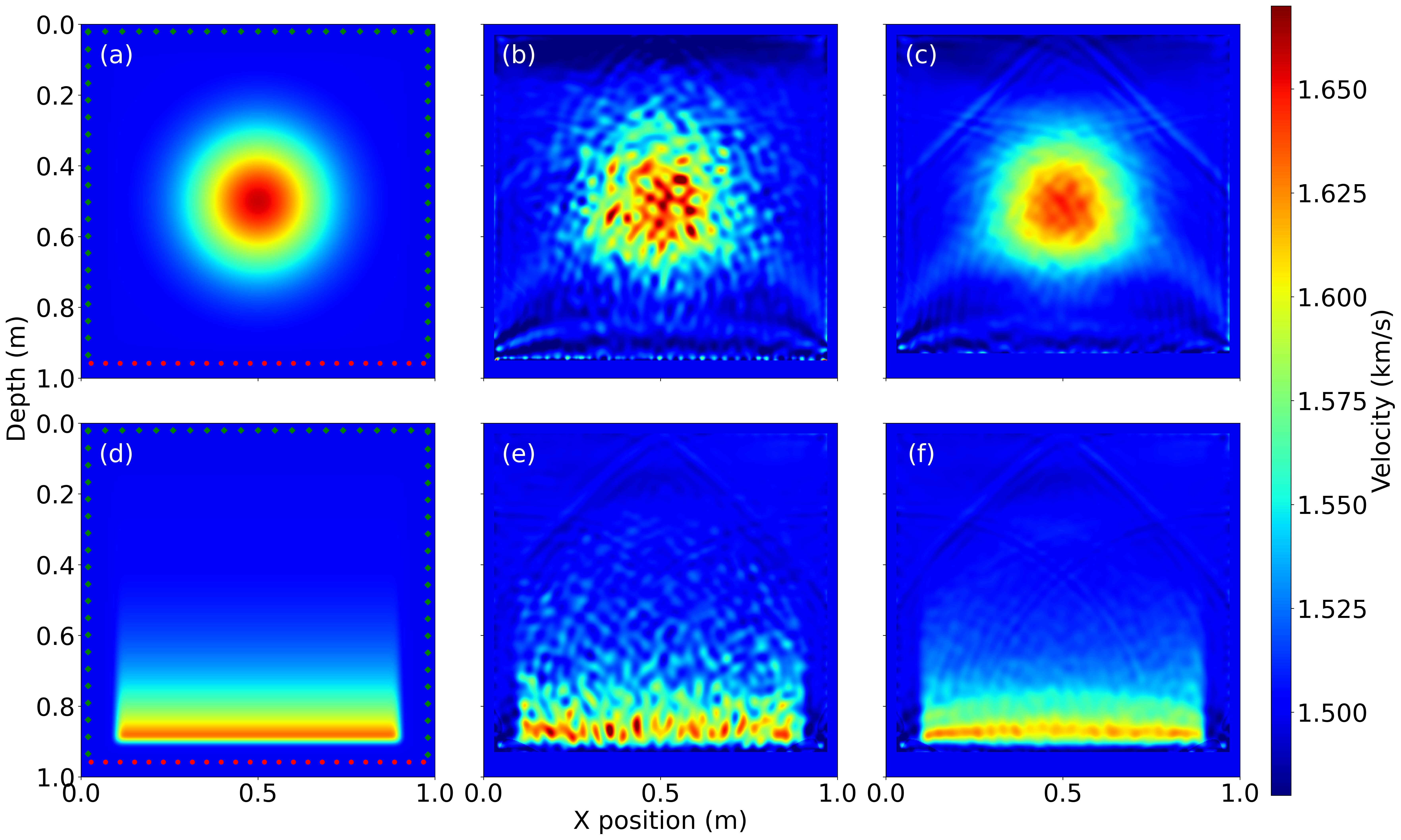}}
    \vspace*{8pt}
    \caption{The true model, and inverted models with and without shot averaging technique based on FWI with $W^2$ objective function. The top row shows the results for the Gaussian function, while the bottom row shows the results for the sediment profile from \cite{chiuMathematicalModelsDistribution2000}. (a) True effective model. (b) Inverted model without shot averaging. (c) The inverted model with shot averaging. (d) True effective model. (e) Inverted model without shot averaging. (f) The inverted model with shot averaging.}
    \label{fig_inverse_het}
\end{figure}

The true model and the inverted model are shown in \cref{fig_inverse_het}. The top row in \cref{fig_inverse_het} shows the results for the Gaussian function, while the bottom row shows the results for the sediment profile from \cite{chiuMathematicalModelsDistribution2000}. The left column shows the true effective model, while the middle column shows the inverted model without the shot averaging technique. The right column shows the inverted model with the shot averaging technique.

The inverted model without shot averaging exhibits oscillatory artifacts due to the high-frequency nature of the data from the heterogeneous medium, which complicates the inversion process. These oscillations contribute to noise in the inverted model. However, as shown in the right column of \cref{fig_inverse_het}, the shot averaging technique effectively reduces these oscillatory noises by averaging out random fluctuations, improving the stability and reliability of the inversion. With shot averaging, the inverted model provides a more accurate approximation of the true model, highlighting the technique's effectiveness in mitigating oscillatory artifacts.

\subsection{Inversion with model mollification}
Now we consider the inversion with model mollification. The mollified kernel is chosen as the Gaussian function with a standard deviation $\sigma_\text{moll}=10,20$ for two cases, respectively.

\begin{figure}[t]
    \centerline{\includegraphics[width=0.8\textwidth]{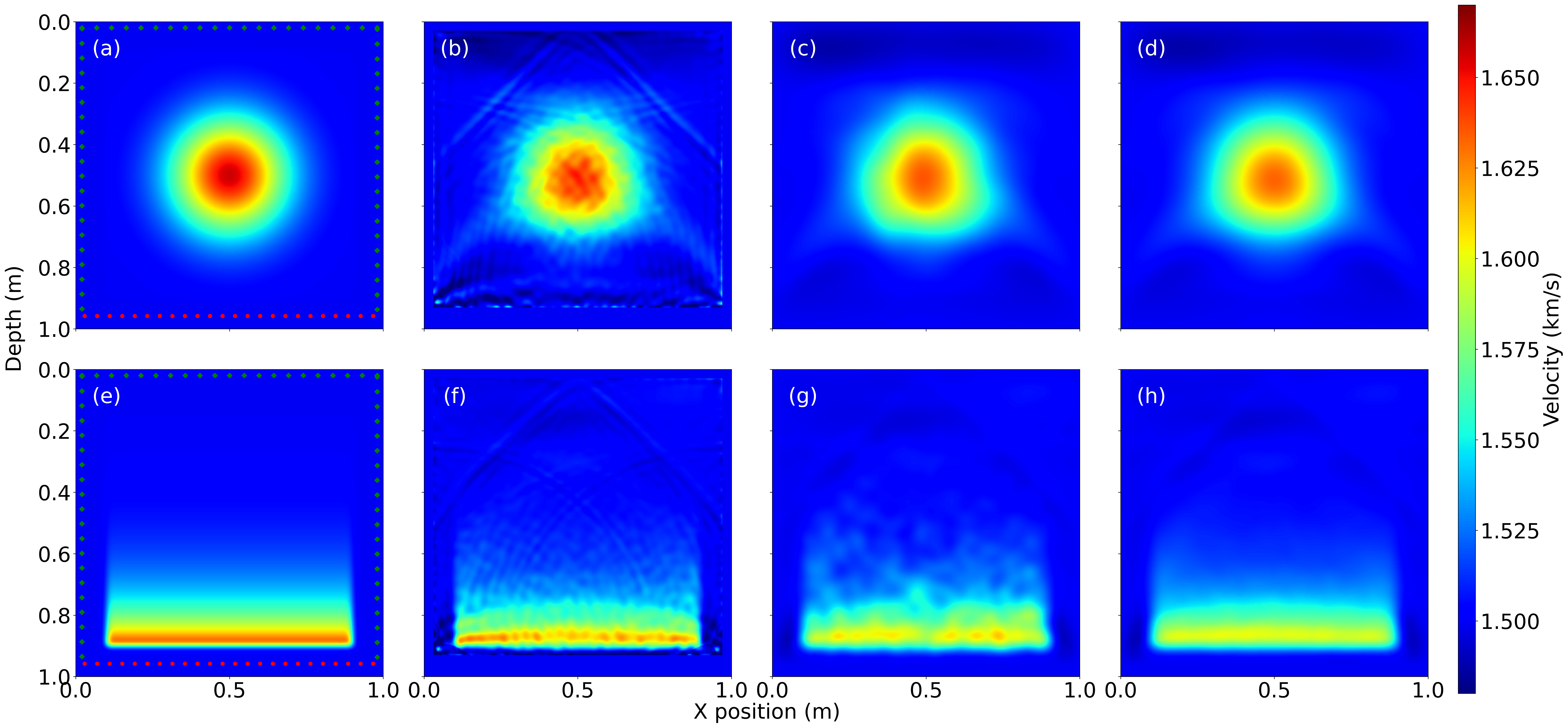}}
    \vspace*{8pt}
    \caption{The true model, and inverted models with model mollification based on FWI with $W^2$ objective function. The top row shows the results for the Gaussian function, while the bottom row shows the results for the sediment profile from \cite{chiuMathematicalModelsDistribution2000}. (a) True effective model. (b) Inverted model with shot averaging. (c) The inverted model with model mollification. (d) The inverted model with shot averaging and model mollification. (e) True effective model. (f) Inverted model with shot averaging. (g) The inverted model with model mollification. (h) The inverted model with shot averaging and model mollification.}
    \label{fig_inverse_moll}
\end{figure}

The inverted results are shown in \cref{fig_inverse_moll}. The top row shows the results for the Gaussian function, while the bottom row shows the results for the sediment profile from \cite{chiuMathematicalModelsDistribution2000}. The first column shows the true effective model, the second column shows the inverted model with shot averaging, the third column shows the inverted model with model mollification, and the last column shows the inverted model with both shot averaging and model mollification.

One can observe that the oscillatory artifacts in the inverted model are effectively mitigated by the model mollification. The mollified model exhibits a smoother profile, which is more consistent with the true model. Even inversion with only one measurement, the model mollification can still provide a good reconstruction of the true model. This demonstrates the effectiveness of the approach of the model mollification in improving the accuracy of the inversion results. In practical terms, this approach can significantly reduce the computational cost associated with data acquisition, as it eliminates the need for multiple realizations of the heterogeneous medium. Furthermore, the model mollification technique can be easily integrated into existing inversion frameworks, thereby enhancing the precision and reliability of the inversion process. Overall, the combination of the shot averaging technique and model mollification can further enhance the stability and accuracy of the inversion results, as demonstrated in the last column of \cref{fig_inverse_moll}.

\subsection{Estimation of the sediment concentration}
The results obtained from the sediment concentration estimation are presented in \cref{tab:sediment_concentration_gaussian,tab:sediment_concentration_chiu}. Here we list 
\begin{enumerate}
    \item the true sediment concentrations:
    \begin{enumerate}
        \item the exact sediment concentration calculated from the single sediment sample,
        \item the estimated sediment concentration from the probability pattern $p(x)$ which can be seen as the target value of the inversions,
        \item the exact sediment concentration calculated from multiple sediment samples,
    \end{enumerate}
    \item the measured sediment concentration from the inverted medium model $\tilde{p}(x)$ using four different methods,
    \item the relative error between the true value and the measured value.
\end{enumerate}

The reason why we compare the inverted sediment concentrations with three target values is that the true sediment distribution is obtained from a single realization of the probability pattern $p(x)$, which is not quite reliable. This sediment concentration would vary from one realization to another, while the estimated sediment concentration from the probability pattern $p(x)$ and the average sediment concentration is the mean values of the sediment concentration of all the realizations, which are more reliable. Hence, the comparison with the estimated values from the probability pattern $p(x)$ and the exact values calculated from the multiple realizations can provide a more comprehensive evaluation of the inversion results.

The results demonstrate a strong agreement between the estimated and actual sediment concentrations, highlighting the accuracy of the estimation process. Furthermore, the application of the model mollification method significantly improves the accuracy of the sediment concentration estimates, as reflected in the reduced relative errors. However, it is important to note that the relative errors cannot be completely eliminated due to the zeroth order approximation inherent in the effective medium model.

Moreover, it is important to note that the shot averaging technique would not necessarily decrease the relative error of the estimated sediment concentration. As analyzed in \cref{sec_solveIP}, averaging over independent realizations suppresses the realization-dependent random fluctuation, but it does not eliminate the systematic bias at the fixed microscale. This explains why, in our numerical tests, shot averaging alone yields only limited improvement in the reconstruction, especially when the domain of sediment occurrence is small, such as the model from \cite{chiuMathematicalModelsDistribution2000}. Nevertheless, it is still beneficial to apply the shot averaging technique for the case of large domain of sediment occurrence, such as the Gaussian model, as the reduction of random fluctuations becomes noticeable and leads to a modest error decrease.

Another observation is that the combination of the shot averaging technique and the model mollification method does not necessarily provide the best inversion results. We speculate that the combination of these two techniques may lead to over-smoothing of the inverted model, which could potentially result in the loss of some high-frequency information and the local peaks in the inverted model. This could explain the relatively higher relative errors of the estimated sediment concentration in the last row of \cref{tab:sediment_concentration_gaussian,tab:sediment_concentration_chiu}.

Overall, the results demonstrate the effectiveness of the proposed approach for sediment concentration estimation and suggest its potential for practical applications.

\begin{table}[ht]
    \centering
    \renewcommand\arraystretch{1.5}
    \begin{tabular}{|*{5}{c|}}
        \hline
        \multicolumn{2}{|c|}{}& \makecell{Exact \\ (single)} & \makecell{Estimated \\ from $p(x)$} & \makecell{Exact \\ (multiple)} \\
        \cline{3-5}
        \multicolumn{2}{|c|}{\multirow{-2.35}{*}{\diagbox[height=3.5\line,width=6cm]{Methods and inverted values}{True values}}}& 0.0277& 0.0294& 0.0279\\
        \hline
        Without shot averaging   & 0.0252 & 8.7918\% & 14.3122\% & 9.7560\%\\ \hline
        With only shot averaging  & 0.0260 & 5.9385\% & 11.6316\% &6.9329\%\\ \hline
        \makecell{With only model \\ mollification} & 0.0267 & 3.5592\% & 9.3964\% & 4.5788\%\\ \hline
        \makecell{With shot averaging and \\ model mollification} & 0.0263 & 4.9898\% & 10.7404\% & 5.9942\%\\
        \hline
    \end{tabular}
    \vspace*{8pt}
    \caption{The relative errors of estimated sediment concentration of Gaussian model}
    \label{tab:sediment_concentration_gaussian}
\end{table}
\begin{table}[ht]
    \centering
    \renewcommand\arraystretch{1.5}
    \begin{tabular}{|*{5}{c|}}
        \hline
        \multicolumn{2}{|c|}{}& \makecell{Exact \\ (single)} & \makecell{Estimated \\ from $p(x)$} & \makecell{Exact \\ (multiple)} \\
        \cline{3-5}
        \multicolumn{2}{|c|}{\multirow{-2.35}{*}{\diagbox[height=3.5\line,width=6cm]{Methods and inverted values}{True values}}}& 0.0251 & 0.0268 & 0.0255\\
        \hline
        Without shot averaging  & 0.0233 & 7.0658\% & 13.1448\% & 8.4799\%\\ \hline
        With only shot averaging & 0.0230 & 8.3619\% & 14.3562\% &9.7563\%\\ \hline
        \makecell{With only model \\ mollification}  & 0.0240 & 4.2620\% & 10.5245\% &5.7188\%\\\hline
        \makecell{With shot averaging and \\ model mollification} & 0.0231 & 7.7666\% & 13.7998\% &9.1700\% \\
        \hline
    \end{tabular}
    \vspace*{8pt}
    \caption{The relative errors of estimated sediment concentration of model from \cite{chiuMathematicalModelsDistribution2000}}
    \label{tab:sediment_concentration_chiu}
\end{table}

\section{Conclusion}\label{sec_conclusions}
In this work, we have developed a multiscale framework for estimating sediment concentration from acoustic wave measurements. The microscopic sediment configuration is modeled by a spatially inhomogeneous Poisson cloud, which provides a direct connection between the random particle distribution and the macroscopic probability $p(x)$ of sediment occurrence. For the associated heterogeneous acoustic wave equation, we derived an effective model and established a quantitative homogenization estimate. The resulting effective slowness square is explicitly determined by $p(x)$, thereby linking the stochastic microscopic model to the macroscopic quantity relevant for sediment concentration measurement.

Based on this effective description, we formulate the recovery of the sediment distribution as an inverse medium problem from partial boundary measurements. The use of the effective model avoids resolving individual particles in the inversion and therefore substantially reduces the computational cost. Numerical experiments show that the effective model captures the dominant macroscopic features of wave propagation in the heterogeneous medium and that the reconstructed effective coefficient can be used to obtain reasonable estimates of sediment concentration. The experiments also illustrate that model mollification and averaging over microscopic realizations can reduce oscillatory effects arising from the fine-scale random medium.

The present framework uses a leading-order effective model and focuses on the spatially inhomogeneous Poisson-cloud description of the sediment particles. Several extensions remain of interest. Higher-order effective models could be used to capture wave-field fluctuations and multiple-scattering effects that are absent from the leading-order approximation. It would also be useful to extend the quantitative analysis to broader classes of spatially varying random media and to investigate systematically how the reconstruction depends on the separation between the microscopic scale, the inversion mesh size, and the acquisition geometry. These questions may further clarify the range of validity and practical performance of the effective-model approach for sediment concentration measurement.

\bibliographystyle{abbrv}
\bibliography{refs} 

\end{document}